\documentclass[10pt]{amsart}

\usepackage{anysize}
\usepackage[utf8x]{inputenc}
\usepackage[english]{babel}
\usepackage{amssymb, amsmath, amsthm}
\usepackage{mathabx}
\usepackage{mathtools}
\usepackage{tikz}
\usepackage{tikz-cd}
\usetikzlibrary{arrows}
\usepackage{graphicx}
\usepackage{url}
\usepackage{color}
\usepackage[all,cmtip]{xy}
\usepackage[hidelinks,bookmarksnumbered]{hyperref}
\usepackage{cleveref}
\usepackage{enumitem}

\let\cal\mathcal
\def\AA{{\cal A}}
\def\BB{{\cal B}}
\def\CC{{\cal C}}
\def\DD{{\cal D}}
\def\EE{{\cal E}}
\def\FF{{\cal F}}

\let\blb\mathbb
 
\def\bC{{\blb C}} 
\def\bD{{\blb D}}
\def\bE{{\blb E}}
\def\bI{{\blb I}}

\def\bF{{\blb F}}

\title{On the obscure axiom for one-sided exact categories}
\author{Ruben Henrard}
\address{Ruben Henrard \\ Universiteit Hasselt \\ Campus Diepenbeek \\ Departement WNI \\ 3590 Diepenbeek \\ Belgium}
\email{ruben.henrard@uhasselt.be}
\author{Adam-Christiaan van Roosmalen}
\address{Adam-Christiaan van Roosmalen \\ Universiteit Hasselt \\ Campus Diepenbeek \\ Departement WNI \\ 3590 Diepenbeek \\ Belgium}
\email{adamchristiaan.vanroosmalen@uhasselt.be}

\newtheorem{theorem}{Theorem}[section]
\newtheorem{proposition}[theorem]{Proposition}
\newtheorem{lemma}[theorem]{Lemma}
\newtheorem{corollary}[theorem]{Corollary}

\theoremstyle{definition}
\newtheorem{definition}[theorem]{Definition}
\newtheorem{remark}[theorem]{Remark}
\newtheorem{notation}[theorem]{Notation}
\newtheorem{example}[theorem]{Example}

\DeclareMathOperator{\inflation}{\rightarrowtail}
\DeclareMathOperator{\deflation}{\twoheadrightarrow}

\DeclareMathOperator{\im}{im}
\DeclareMathOperator{\coim}{coim}
\DeclareMathOperator{\coker}{coker}

\DeclareMathOperator{\Hom}{Hom}
\DeclareMathOperator{\Homex}{Hom_{\mathrm{exact}}}
\DeclareMathOperator{\Mor}{Mor}

\DeclareMathOperator{\modd}{mod}

\DeclareMathOperator{\rep}{rep}

\DeclareMathOperator{\Ab}{\mathsf{Ab}}

\DeclareMathOperator{\Ob}{Ob}

\DeclareMathOperator{\Ac}{\textbf{Ac}^b}

\DeclareMathOperator{\Cb}{\textbf{C}^b}
\DeclareMathOperator{\Db}{\textbf{D}^b}
\DeclareMathOperator{\Kb}{\textbf{K}^b}

\newcommand{\Rt}{\textbf{R3}}
\newcommand{\Rtp}{{\textbf{R3}^\textbf{+}}}
\newcommand{\Rtm}{{\textbf{R3}^\textbf{-}}}
\newcommand{\Lt}{\textbf{L3}}
\newcommand{\Ltp}{{\textbf{L3}^\textbf{+}}}

\newcommand{\bDm}{\bD_{\textrm{max}}}
\newcommand{\bIm}{\bI_{\textrm{max}}}

\newcommand{\ex}[1]{#1^{\textrm{ex}}}

\newcommand{\fD}{\mathfrak{D}}
\newcommand{\sfD}{\mathfrak{sD}}
\newcommand{\fDp}{\mathfrak{D}^+}
\newcommand{\fDm}{\mathfrak{D}^-}

\newcommand{\fI}{\mathfrak{I}}

\newcommand{\fE}{\mathfrak{E}}

\makeatletter
\newcommand{\myitem}[1]{%
\item[#1]\protected@edef\@currentlabel{#1}%
}
\makeatother


\begin{document}

\begin{abstract}
	One-sided exact categories are obtained via a weakening of a Quillen exact category.  Such one-sided exact categories are homologically similar to Quillen exact categories: a one-sided exact category $\EE$ can be (essentially uniquely) embedded into its exact hull $\ex{\EE};$ this embedding induces a derived equivalence $\Db(\EE) \to \Db(\ex{\EE})$.
	
	Whereas it is well known that Quillen's obscure axioms are redundant for exact categories, some one-sided exact categories are known to not satisfy the corresponding obscure axiom.  In fact, we show that the failure of the obscure axiom is controlled by the embedding of $\EE$ into its exact hull $\ex{\EE}.$
	
	In this paper, we introduce three versions of the obscure axiom (these versions coincide when the category is weakly idempotent complete) and establish equivalent homological properties, such as the snake lemma and the nine lemma.  We show that a one-sided exact category admits a closure under each of these obscure axioms, each of which preserves the bounded derived category up to triangle equivalence.
\end{abstract}

\maketitle

\tableofcontents

\section{Introduction}

Quillen introduced the notion of an exact category in \cite{Quillen73} as a framework for homological algebra and algebraic $K$-theory. An exact category is an additive category together with a chosen class of kernel-cokernel pairs (called \emph{conflations}) satisfying 8 axioms. The kernel morphism of a conflation is called an \emph{inflation} and the cokernel morphism is called a \emph{deflation}. The 8 axioms can be partitioned into two dual sets of axioms, with axioms \ref{R0}-\ref{R3} referring solely to the deflation side and axioms \ref{L0}-\ref{L3} referring solely to the inflation side (see definition \ref{definition:OneSidedExactCategory} for the axioms).  It is well known that the axioms of an exact category are not minimal; in particular, the \emph{obscure axioms} (as they were referred to in \cite{ThomasonTrobaugh90} and subsequently in \cite{Buhler10}) \ref{L3} and \ref{R3} are superfluous (see \cite{Keller90,Yoneda60}).

Requiring only the axioms on the deflation side (axioms \ref{R0}-\ref{R2}) gives the definition of a \emph{deflation-exact} category.  These axioms imply that the deflations define a Grothendieck pretopology on the category (see \cite{Rosenberg11} or \cite{KaledinLowen15}).  Requiring the inflation side of the axioms gives rise to an inflation-exact category.  Such one-sided exact categories are relative versions of one-sided quasi-abelian  categories (also called one-sided almost abelian categories), and have been introduced by \cite{BazzoniCrivei13} and \cite{Rump11} (in \cite{Rump11}, the obscure axiom is part of the definition).

Recently, the theory of one-sided exact categories has gained some interest.  One-sided exact categories appear as an intermediate step in the construction of the maximal exact structure on an additive category (\cite{Crivei12, RichmanWalker77, Rump11, SiegWegner11}); similarly, they appear as an intermediate step in a recent approach to the $K$-theory of the category of locally compact abelian groups (\cite{BraunlingHenrardvanRoosmalen20}, shortening the proof from \cite{Braunling20}, see also \cite{Clausen17}), as a localization of an exact category (\cite{HenrardVanRoosmalen19b, HenrardVanRoosmalen19a}).  Furthermore, natural examples of one-sided exact categories can be found in representation theory (\cite{CaenepeelVanOystaeyen19, HenrardvanRoosmalen20}), functional analysis (for example, the category of $\mathsf{LB}$-spaces or the category of complete Hausdorff locally convex spaces, see \cite{HassounShahWegner20}).  An additive homological category \cite{BorceuxBourn04} is an example of a one-sided exact category, as is a so-called one-morphism Grothendieck pretopology (\cite{KaledinLowen15, Rosenberg11}).

In this paper, we consider an additive category $\AA$ together with a class of conflations $\bD$ such that the pair $\EE = (\AA, \bD)$ is a one-sided exact category.  We will state our results for deflation-exact categories, leaving the dual statements for inflation-exact categories to the reader.  In this introduction, and throughout most of the paper, we assume that $\bD$ contains all split kernel-cokernel pairs (this condition has been referred to in \cite{BazzoniCrivei13} as axiom \ref{R0*}, and we will keep this terminology).

At first glace, a one-sided exact category might seem a considerably weaker structure than an exact category.  Closer inspection reveals that a one-sided exact category $\EE$ admits an exact hull $\overline{\EE}$, meaning that there exists an exact and full embedding $j\colon \EE \to \ex{\EE}$ of $\EE$ into an exact category $\ex{\EE}$, which is 2-universal among all exact functors from $\EE$ to an exact category (see \cite{HenrardVanRoosmalen19b, Rosenberg11}).  Moreover, the embedding $j$ lifts to a derived equivalence $\Db(\EE) \to \Db(\ex{\EE}).$  Hence, a one-sided exact structure is (homologically and $K$-theoretically) close to a Quillen exact structure.  However, to translate homological properties from the exact hull $\ex{\EE}$ to $\EE$, the following property would be useful: a sequence $X \to Y \to Z$ is a conflation in $\EE$ if and only if $j(X) \to j(Y) \to j(Z)$ is a conflation in $\ex{\EE}.$  In this paper, we show that this property is equivalent to the obscure axiom (see theorem \ref{theorem:EquivalentFormulations} or proposition \ref{proposition:WhenConflationInCompletionIntroduction}); this illustrates how the obscure axiom implies the existence of homological properties close to those of an exact category.

However, some one-sided exact categories of interest do not satisfy the corresponding obscure axiom.  For example, the conflation structure given by all semi-stable cokernels \cite{Rump15} and the quotient of an exact category by a deflation-percolating subcategory \cite{HenrardVanRoosmalen19a} are examples of deflation-exact categories, possibly not satisfying the obscure axiom.

\subsection{Homological consequences of the obscure axiom}

Let $\EE = (\AA, \bD)$ is a deflation-exact category.  We say that $\EE$ satisfies the obscure axiom \ref{R3} if the following condition holds: for any morphism $f\colon X \to Y$, if $f$ has a kernel and there is a morphism $g\colon X' \to X$ such that $g \circ f \colon  X' \to X \to Y$ is a deflation, then $f$ is a deflation.  In addition, we also consider the following two versions of the obscure axiom by altering the condition ``$f$ has a kernel'': if we instead require $f$ to admit all pullbacks, we obtain axiom \ref{R3-}; if we remove the condition on the kernel of $f$ altogether, we obtain axiom \ref{R3+}.

Following \cite{BazzoniCrivei13}, we call a deflation-exact category \emph{strong} if it satisfies the obscure axiom \ref{R3}.  The following equivalent conditions illustrate that a strongly deflation-exact category is homologically close to an exact category..

\begin{theorem}\label{theorem:EquivalentFormulations}
	Let $\EE$ be a deflation-exact category satisfying axiom \ref{R0*}. The following are equivalent:
	\begin{enumerate}
		\item\label{item:EquivalentFormulationsA} Axiom \ref{R3} holds.
		\item\label{item:EquivalentFormulationsB} The nine lemma holds (see theorem \ref{theorem:NineLemma} for a precise statement of the nine lemma).
		\item\label{item:EquivalentFormulationsE} Conflations are closed under retracts.
		\item\label{item:EquivalentFormulationsF} Conflations are closed under direct summands.
		\item\label{item:EquivalentFormulationsC} If $\begin{psmallmatrix}0&g\end{psmallmatrix}\colon X\oplus Y\deflation Z$ is a deflation and $g$ admits a kernel, then $g$ is a deflation.
		\item\label{item:EquivalentFormulationsD} A morphism $g\colon Y \to Z$ with kernel $f\colon X\to Z$ is a deflation if and only if there exists a deflation $f'\colon Y' \deflation Y$ such that $g \circ f$ is a deflation.
		\item\label{item:EquivalentFormulationsH} A sequence $X\xrightarrow{f}Y\xrightarrow{g}Z$ is a conflation in $\EE$ if and only if there is a triangle $i(X)\xrightarrow{i(f)}i(Y)\xrightarrow{i(g)}i(Z)\rightarrow \Sigma i(X)$ in $\Db(\EE)$.
		\item\label{item:EquivalentFormulationsI} A sequence $X\xrightarrow{f}Y\xrightarrow{g}Z$ is a conflation in $\EE$ if and only if it is a conflation in the exact hull $\overline{\EE}$ (see theorem \ref{Theorem:ExactHullIntroduction}).
		\item\label{item:EquivalentFormulationsJ} The embedding $\EE \to \widehat{\EE}$ of $\EE$ into its weak idempotent completion $\widehat{\EE}$ is fully exact.
	\end{enumerate}
\end{theorem}

For a deflation-exact category $\EE$, axiom \ref{R3+} is equivalent to satisfying the obscure axiom \ref{R3} and being \emph{weakly idempotent complete} (i.e.~every retract has a kernel; this property is also called \emph{divisive}, see \cite{Rump20}): deflation-exact categories satisfying axiom \ref{R3+} are closely related to a weakly idempotent complete exact category.  Specifically, we have the following theorem. 

\begin{theorem}\label{theorem:WeaklyIdempotentCompleteEquivalentFormulations}
	Let $\EE$ be a deflation-exact category satisfying axiom \ref{R0*}. The following are equivalent:
	\begin{enumerate}
		\item\label{item:WeaklyIdempotentCompleteEquivalentFormulationsA} Axiom \ref{R3+} holds.
		\item\label{item:WeaklyIdempotentCompleteEquivalentFormulationsB} $\EE$ is weakly idempotent complete and axiom \ref{R3} holds.
		\item\label{item:WeaklyIdempotentCompleteEquivalentFormulationsC} If $\begin{psmallmatrix}0&g\end{psmallmatrix}\colon X\oplus Y\deflation Z$ is a deflation, then $g$ is a deflation.
		\item\label{item:WeaklyIdempotentCompleteEquivalentFormulationsF} A morphism $g\colon Y \to Z$ is a deflation if and only if there exists a deflation $f\colon X \deflation Y$ such that $g \circ f$ is a deflation.
		\item\label{item:WeaklyIdempotentCompleteEquivalentFormulationsD} The ker-coker-sequence property holds (see proposition \ref{proposition:Ker-Coker-Sequence}).
		\item\label{item:WeaklyIdempotentCompleteEquivalentFormulationsE} The (short) snake lemma holds (see corollary \ref{corollary:SnakeLemmaI} and theorem \ref{theorem:SnakeLemmaII}).
				\item\label{item:WeaklyIdempotentCompleteEquivalentFormulationsH} Deflations are closed under retracts, i.e.~for any commutative diagram
		\[\xymatrix{
		Y \ar@{->>}[r]^{g} \ar[d] & Z \ar[d] \\
		{Y'} \ar[r]^{g'} & {Z'}		}\]
		where the vertical arrows are retractions, if $g$ is a deflation, then $g'$ is a deflation.
		\item\label{item:WeaklyIdempotentCompleteEquivalentFormulationsG} Deflations are closed under direct summands, i.e.~for any two morphisms $g\colon Y \to Z$ and $g'\colon Y'\to Z'$, if the direct sum $Y \oplus Y'\deflation Z \oplus Z'$ is a deflation, then $g$ and $g'$ are deflations.%
	\end{enumerate}
\end{theorem}

Both theorem \ref{theorem:EquivalentFormulations} and theorem \ref{theorem:WeaklyIdempotentCompleteEquivalentFormulations} will be proved in \S\ref{section:EquivalentFormulations}.

\subsection{The obscure closure of a one-sided exact category}

Let $\EE = (\AA, \bD)$ be a deflation-exact category (such that the split kernel-cokernel pairs are conflations).  It was shown in \cite{Rump11} that $\bD$ can be restricted to a largest strongly deflation-exact substructure. This restriction usually alters the derived category and might as such be undesirable.   In the opposite direction, it need not be possible to extend $\bD$ to a strongly deflation-exact structure on $\AA$.  The obstruction lies with conflations of the form
\[\xymatrix@1{{X \oplus A\, } \ar@{>->}[r]^{\begin{psmallmatrix}i & 0 \\ 0 & 1\end{psmallmatrix}} & Y \oplus A \ar@{->>}[r]^-{\begin{psmallmatrix} p & 0 \end{psmallmatrix}} & Z};\]
even though the pullback $P$ of a morphism $f\colon Z' \to Z$ along $\begin{psmallmatrix} p & 0 \end{psmallmatrix}\colon  Y \oplus A \deflation Z$ may exist, the pullback of $f$ along $p\colon Y  \to Z$ need not exist in $\AA$. 

This leaves the following options: add only those conflations for which the required pullbacks exist, or change the underlying category $\AA$.  The former approach yields the closure under axiom \ref{R3-}, while the latter approach yields the closure under axiom \ref{R3} or \ref{R3+}.

\begin{theorem}\label{theorem:IntroductionObscureClosures}
	Let $\EE = (\AA, \bD)$ be a deflation-exact category such that the split kernel-cokernel pairs are conflations.
	\begin{enumerate}
		\item There is a smallest deflation-exact structure $\bD_{\Rtm} \supseteq \bD$ on $\AA$ satisfying axiom \ref{R3-}.
		\item There is an exact functor $\EE \to \EE_\Rt$ from $\EE$ to a deflation-exact category $\EE_\Rt$ satisfying axiom \ref{R3}, universal among all exact functors to deflation-exact categories satisfying axiom \ref{R3}.
		\item There is an exact functor $\EE \to \EE_\Rtp$ from $\EE$ to a deflation-exact category $\EE_\Rtp$ satisfying axiom \ref{R3+}, universal among all exact functors to deflation-exact categories satisfying axiom \ref{R3+}.
	\end{enumerate}
	Moreover, the embeddings of $\EE$ into each of these three closures are derived equivalences.
\end{theorem}

The deflation structure $\EE_\Rtm = (\AA, \bD_\Rtm)$ is the largest extension of $\bD$ such that $\EE \to \EE_\Rtm$ is a derived equivalence.  The categories $\EE_\Rt$ and $\EE_\Rtp$ need not have $\AA$ as underlying category: the underlying additive category of $\EE_\Rtp$ is the weak idempotent completion $\widehat{\AA}$ of $\AA$, and the underlying additive category of $\EE_\Rt$ is a subcategory of $\widehat{\AA}$ containing $\AA$.

As the obscure axiom yields a slew of interesting homological properties, it is useful to know which sequences in $\EE$ become conflations in $\EE_\Rt$ (or $\EE_\Rtp$ or $\ex\EE$).  The following proposition addresses this (see proposition \ref{proposition:WhenConflationInCompletion} in the text).

\begin{proposition}\label{proposition:WhenConflationInCompletionIntroduction}
Let $X\xrightarrow{f} Y \xrightarrow{g} Z$ be a sequence in a deflation-exact category $\EE$.  The following are equivalent:
\begin{enumerate}
	\item the sequence is a conflation in $\EE_\Rt$,
	\item the sequence is a conflation in $\EE_\Rtp$,
	\item the sequence is a conflation in $\ex{\EE}$,
	\item there is an $A \in \EE$ such that $\xymatrix@1{ {X\oplus A\,}\ar@{>->}[r]^{\begin{psmallmatrix} f & 0 \\ 0 & 1 \end{psmallmatrix}} & {Y \oplus A} \ar@{->>}[r]^-{\begin{psmallmatrix} g & 0 \end{psmallmatrix}} & Z}$ is a conflation in $\EE.$
\end{enumerate}
\end{proposition}  

In \S\ref{Section:Lattices}, we apply these insights to compare the different lattice of one-sided exact structures on an additive category $\AA$ and its weak idempotent completion $\widehat{\AA}$.  Specifically, we show that there is a Galois connection between the deflation-exact structures satisfying axiom \ref{R3-} on $\AA$ and the deflation-exact structures satisfying axiom \ref{R3+} on $\widehat{\AA}.$

\subsection{The exact structure of stable kernel-cokernel pairs}

One-sided exact structures have been used to construct the maximal exact structure on an additive category $\AA$: for a weakly idempotent complete category, the maximal exact structure is given by intersecting the maximal deflation-exact structure $\bDm$ and maximal inflation-exact structure $\bIm$ on $\AA$ (see \cite{RichmanWalker77,SiegWegner11} for the pre-abelian case and \cite{Crivei12} for the more general weakly idempotent complete case).  The conflations in $\bDm \cap \bIm$ are the stable kernel-cokernels.

When $\AA$ is not weakly idempotent complete, the stable kernel-cokernels need not form an exact structure.  Here, the maximal exact structure is obtained by intersecting the largest strongly deflation- and strongly inflation-exact categories, but this set of conflations might be considerably smaller than the set of stable kernel-cokernels (see \cite{Rump11, Rump15}).

In \S\ref{Section:StableConflations}, we show that, given any additive category $\AA$, the intersection $\bI \cap \bD$ of an inflation-exact structure $\bI$ and a deflation-exact structure $\bD$ ``almost'' endows $\AA$ with the structure of an exact category: even though $(\AA, \bI \cap \bD)$ is not an exact category, the weak idempotent completion $(\widehat{\AA}, \widehat{\bI \cap \bD})$ is an exact category.  This suggests that an exact closure of $(\AA, \bI \cap \bD)$ can be constructed inside $(\widehat{\AA}, \widehat{\bI \cap \bD})$.  The following theorem is proposition \ref{proposition:UniversalPropertyIntersection}.

\begin{theorem}
Let $\bD$ and $\bI$ be a deflation- and an inflation-exact structure on an additive category $\AA.$  There is an exact category $(\AA', \bE)$ and a conflation-exact functor $\varphi\colon (\AA, \bD \cap \bI) \to (\AA', \bE')$ satisfying the following 2-universal property: for each exact category $(\BB, \EE)$, the functor
\[-\circ \varphi\colon \Homex((\AA,\bD \cap \bI),(\BB,\bF)) \to \Homex((\AA',\bE'), (\BB,\bF))\]
is an equivalence. 
\end{theorem}

The additive category $\AA'$ in this theorem is found as a full subcategory of the weak idempotent completion $\widehat{\AA}$, and contains the category $\AA$.  Put differently, $\AA \subseteq \AA' \subseteq \widehat{\AA}.$

\subsection*{Acknowledgments} The second author is currently a postdoctoral researcher at FWO (12.M33.16N).
\section{Preliminaries}\label{section:Preliminaries}

This section is preliminary in nature.  We recall the definition of a one-sided exact category as given in \cite{BazzoniCrivei13}.  Note that our left exact categories are called right exact in \cite{BazzoniCrivei13} and vice-versa. To avoid confusion, we will refer to one-sided exact categories by either inflation-exact categories or deflation-exact categories referring directly to the underlying axioms instead.

\subsection{Characterization of pullbacks and pushouts}

We start with some general statements about pullbacks and pushouts.  The following statement is \cite[proposition~I.13.2]{Mitchell65} together with its dual.

\begin{proposition}\label{proposition:MitchellPushout}\label{proposition:MitchellPullback}
Let $\CC$ be any pointed category.
\begin{enumerate}
  \item\label{enumerate:MitchellPullback} Consider a  diagram
	\[ \xymatrix{ {X'} \ar[r]^{f'} & {Y'} \ar[d]^{h} & \\ {X} \ar[r]^{f} & {Y} \ar[r]^{g} & Z } \]
where $f$ is the kernel of $g$.  The left-hand side can be completed to a pullback square if and only if $f'$ is the kernel of $gh$.
\item\label{enumerate:MitchellPushout} Consider a diagram
	\[ \xymatrix{ {X} \ar[r]^{f} & {Y} \ar[d]^{h} \ar[r]^{g} & Z \\ & {Y'} \ar[r]^{g'} & Z'} \]
where $g$ is the cokernel of $f$.  The right-hand side can be completed to a pushout square if and only if $g'$ is the cokernel of $hf$.
\end{enumerate}
\end{proposition}

The following well-known proposition states that pullbacks preserve kernels and pushouts preserve cokernels (see \cite[lemma 1]{Rump11}).

\begin{proposition}\label{proposition:PullbacksPreserveKernels}
Let $\CC$ be any pointed category.  Consider the following commutative diagram in $\CC$:
\[\xymatrix{
A \ar[d]^{g} \ar[r]^{f} & B \ar[d]^{h} \\
C \ar[r]^{f'} & D
}\]
\begin{enumerate}
	\item Assume that the commutative square is a pullback.  The morphism $f$ admits a kernel if and only if $f'$ admits a kernel.  In this case, the composition $\ker(f) \to A \to C$ is the kernel of $f'.$
  \item Assume that the commutative square is a pushout.  The morphism $f$ admits a cokernel if and only if $f'$ admits a cokernel.  In this case, the composition $B \to D \to \coker(f')$ is the cokernel of $f.$
\end{enumerate}
\end{proposition}

\subsection{One-sided exact categories}

\begin{definition}
  A \emph{conflation structure} on an additive category $\AA$ is a chosen class $\bC$ of kernel-cokernel pairs, called \emph{conflations}, such that this class is closed under isomorphisms. The kernel part of a conflation is called an \emph{inflation} and the cokernel part of a conflation is called a \emph{deflation}.  We depict inflations by the symbol $\inflation$ and deflations by $\deflation$.
	
	A \emph{conflation category} is a pair $(\AA, \bC)$ where $\bC$ is a conflation structure on the additive category $\AA.$  A functor $F\colon (\AA_1, \bC_1) \to (\AA_2, \bC_2)$ between conflation categories is called \emph{exact} or \emph{conflation-exact} if it maps conflations to conflations.
\end{definition}

For a conflation category $\EE = (\AA, \bC)$, we often write $X \in \EE$ for $X \in \AA$.

\begin{definition}\label{definition:OneSidedExactCategory}
	A conflation category $\EE = (\AA, \bC)$ is called \emph{deflation-exact} or \emph{right exact} if $\EE$ satisfies the following axioms:
	\begin{enumerate}[label=\textbf{R\arabic*},start=0]
		\item\label{R0} The identity morphism $1_0\colon 0\rightarrow 0$ is a deflation.
		\item\label{R1} Deflations are closed under composition.
		\item\label{R2} The pullback of any morphism along a deflation exists. Moreover, deflations are stable under pullbacks.
	\end{enumerate}
	Dually, a conflation category $\EE$ is called \emph{inflation-exact} or \emph{left exact} if $\EE$ satisfies the following axioms:
	\begin{enumerate}[label=\textbf{L\arabic*},start=0]
		\item\label{L0} The identity morphism $1_0\colon 0\rightarrow 0$ is an inflation.
		\item\label{L1} Inflations are closed under composition.
		\item\label{L2} The pushout of any morphism along an inflation exists. Moreover, inflations are stable under pushouts.
	\end{enumerate}
\end{definition} 

\begin{definition}\label{definition:ObscureAxioms}
	In addition to the axioms listed above, we discuss the following axioms as well.
	\begin{enumerate}[align=left]
		\myitem{\mbox{\textbf{R0}$^\ast$}}\label{R0*} For any $A\in \EE$, $A\rightarrow 0$ is a deflation.
		\myitem{\mbox{\textbf{R3}}}\label{R3} \hspace{0.175cm}If $i\colon A\rightarrow B$ and $p\colon B\rightarrow C$ are morphisms in $\EE$ such that $p$ has a kernel and $pi$ is a deflation, then $p$ is a deflation.
	\end{enumerate}
	The following axioms are dual.
	\begin{enumerate}[align=left]
		\myitem{\mbox{\textbf{L0}$^\ast$}}\label{L0*} For any $A\in \EE$, $0\rightarrow A$ is an inflation.
		\myitem{\mbox{\textbf{L3}}}\label{L3} \hspace{0.175cm}If $i\colon A\rightarrow B$ and $p\colon B\rightarrow C$ are morphisms in $\EE$ such that $i$ has a cokernel and $pi$ is an inflation, then $i$ is an inflation.
	\end{enumerate}
Following \cite{BazzoniCrivei13}, a deflation-exact category satisfying axiom \ref{R3} is called a \emph{strongly deflation-exact category} or a \emph{strongly right exact category}.  Dually, an inflation-exact category satisfying axiom \ref{L3} is called a \emph{strongly inflation-exact category} or a \emph{strongly left exact category}.
\end{definition}

\begin{remark}\label{remark:Definitions}
	\begin{enumerate}
		\item A Quillen exact category is a conflation category satisfying the axioms in \cref{definition:OneSidedExactCategory,definition:ObscureAxioms}. Yoneda showed that axioms \ref{R3} and \ref{L3} are redundant for exact categories (see \cite[p.~525]{Yoneda60}). In fact, Keller showed that axioms \ref{R0},\ref{R1}, \ref{R2} and \ref{L2} imply the other axioms of an exact category (see \cite[appendix~A]{Keller90}).  Dually, an inflation-exact category satisfying axiom \ref{R2} is an exact category.
		\item For a deflation-exact category, axiom \ref{R0*} translates to all split kernel-cokernel pairs being conflations. Similarly, for inflation-exact categories, axiom \ref{L0*} yields that all split kernel-cokernel pairs are conflations (see \cref{proposition:AxiomR0*} below). In this text, we are solely interested in one-sided exact categories having all split kernel-cokernel pairs as conflations. 
	\end{enumerate}
\end{remark}

We recall some basic properties of inflations and deflations in a deflation-exact category (see, for example, \cite{BazzoniCrivei13} or \cite{HenrardVanRoosmalen19a}).

\begin{proposition}\label{proposition:BasicProperties}
	Let $\EE$ be a deflation-exact category. Then:
	\begin{enumerate}
		\item Every isomorphism is a deflation.
		\item Every inflation is a monomorphism.  An inflation which is an epimorphism is an isomorphism.
		\item Every deflation is an epimorphism.  A deflation which is a monomorphism is an isomorphism.
		\item The class of conflations is closed under direct sums.
	\end{enumerate}
\end{proposition}

\begin{proposition}\label{proposition:AxiomR0*}
	Let $\EE$ be a deflation-exact category the following are equivalent:
	\begin{enumerate}
		\item Axiom \ref{R0*} holds.
		\item All split kernel-cokernel pairs are conflations.
		\item Retractions with kernels are deflations.
		\item Coretractions with cokernels are inflations.
	\end{enumerate}
\end{proposition}

\subsection{The derived category of a one-sided exact category}

The derived category of a one-sided exact category was introduced in \cite{BazzoniCrivei13}.  We recall some relevant definitions and properties from \cite{BazzoniCrivei13, HenrardVanRoosmalen19b}. 

We write $\Cb(\EE)$ for the category of bounded cochain complexes over $\EE$ and $\Kb(\EE)$ for the homotopy category of $\EE$. It is well-known that $\Kb(\EE)$ has the structure of a triangulated category induced by the strict triangles in $\Cb(\EE)$. In order to define the (bounded) derived category of a one-sided exact category, we need to define acyclic complexes.

\begin{definition}
	Let $\EE$ be a conflation category. A morphism $f\colon X\to Y$ is called \emph{admissible} or \emph{strict} if it admits a \emph{deflation-inflation factorization}:
	\[\xymatrix{
		X\ar@{->>}[rd]\ar[rr]^{f} && Y\\
		& I\ar@{>->}[ru]
	}\]
	We denote an admissible morphism by $\xymatrix{X\ar[r]|{\circ}^f & Y}$.
\end{definition}

\begin{remark}
The deflation-inflation factorization $X \deflation I \inflation Y$ above is unique (up to isomorphism), moreover, $\coim(f)\cong I\cong \im(f)$. Clearly an admissible map has a kernel and cokernel: $\ker(f)$ is the kernel of $f\deflation \coim(f)$ and $\coker(f)$ is the cokernel of the map $\im(f)\inflation Y$.
\end{remark}

\begin{definition}
	Let $\EE$ be a conflation category. A sequence 
	\[\xymatrix{
		\cdots\ar[r] & X_{i-1}\ar[r]^{f_{i-1}} & X_i\ar[r]^{f_{i}} & X_{i+1}\ar[r]^{f_{i+1}} & X_{i+2}\ar[r]^{f_{i+2}} & \cdots
	}\]
	is called \emph{exact} or \emph{acyclic} if each $f_j$ is admissible, i.e.~factors as $\iota_j\circ \rho_j$ with deflation $\rho_j\colon X_j\deflation \im(f_j)$ and inflation $\iota_j\colon \im(f_j)\inflation X_{i+1}$ such that $\iota_j=\ker(f_{j+1})$ and $\rho_j=\coker(f_{j-1})$ for all $j$.
	
	The full subcategory of $\Kb(\EE)$ consisting of complexes which are homotopic to an acyclic complex is denoted by $\Ac(\EE)$.
\end{definition}

The following lemma is \cite[lemma~7.2]{BazzoniCrivei13}.

\begin{lemma}
Let $\EE$ be a one-sided exact category. The category $\Ac(\EE)$ is a triangulated subcategory of $\Kb(\EE)$.
\end{lemma}

Hence, we arrive at the following definition.

\begin{definition}
	Let $\EE$ be a one-sided exact category. The \emph{bounded derived category} $\Db(\EE)$ is the Verdier localization $\Kb(\EE)/\Ac(\EE)$.
\end{definition}

The following theorem (see \cite[theorem~1.1]{HenrardVanRoosmalen19b}) summarizes some basic properties of the derived category.

\begin{theorem}\label{theorem:DerivedCategoryBasicProperties}
	Let $\EE$ be a deflation-exact category. 
	\begin{enumerate}
		\item The natural embedding $i\colon\EE\hookrightarrow \Db(\EE)$ is fully faithful.
		\item For all $X,Y\in \EE$ and $n>0$, we have $\Hom_{\Db(\EE)}(\Sigma^n i(X),i(y))=0$.
	\end{enumerate}
	If $\EE$ satifies axiom \ref{R0*}, then
	\begin{enumerate}[resume]
		\item a conflation $X\inflation Y\deflation Z$ lifts to a triangle $i(X)\to i(Y)\to i(Z)\to \Sigma i(X)$ in $\Db(\EE)$.
	\end{enumerate}
\end{theorem}

For a deflation-exact category $\EE$, we write $\ex{\EE}$ for the extension-closure of $\EE$ in its bounded derived category $\Db(\EE).$  We may endow $\ex{\EE}$ with the structure of an exact category in the following way (see \cite{HenrardVanRoosmalen19b}, based on \cite{Dyer05}): a sequence $X \to Y \to Z$ is a conflation in $\ex{\EE} \subseteq \Db(\EE)$ if and only if there is a triangle $i(X) \to i(Y) \to i(Z) \to \Sigma i(X)$.  We call $\ex{\EE}$ the \emph{exact hull} of $\EE.$  We recall the following theorem.

\begin{theorem}\label{Theorem:ExactHullIntroduction}
Let $\EE$ be a deflation-exact category satisfying axiom {\ref{R0*}}. The embedding $j\colon \EE\to \ex{\EE}$ is an exact embedding which is $2$-universal among exact functors to exact categories.  Moreover,
\begin{enumerate}
  \item the embedding $\ex{\EE} \to \Db(\EE)$ lifts to a triangle equivalence $\Db(\ex{\EE}) \simeq \Db(\EE)$, and
	\item the embedding $\EE \to \ex{\EE}$ lifts to a triangle equivalence $\Db({\EE}) \simeq \Db(\ex{\EE})$.
\end{enumerate}
\end{theorem}

\subsection{Weakly idempotent complete categories}

Let $\AA$ be any category.  A morphism $f\colon X \to Y$ is called a \emph{retraction} (or a \emph{coretraction}) if there is a morphism $g\colon Y \to X$ such that $f \circ g = 1_Y$ (or $g \circ f = 1_X$, respectively).  The following proposition is \cite[lemma~7.1]{Buhler10}.

\begin{proposition}
In an additive category $\AA$, the following are equivalent.
\begin{enumerate}
\item Every retraction has a kernel.
\item Every coretraction has a cokernel.
\end{enumerate}
\end{proposition}

A category satisfying the conditions of the previous proposition is called \emph{weakly idempotent complete}. The following is an immediate corollary to \cref{proposition:AxiomR0*}.

\begin{corollary}\label{corollary:RetractionsAreDeflations}
	Let $\EE$ be a  weakly idempotent complete deflation-exact category.  If $\EE$ satisfies axiom \ref{R0*}, then retractions are deflations and coretractions are inflations.
\end{corollary}

Every additive category $\AA$ has a weak idempotent closure $\widehat{\AA}$. Explicitly, the weak idempotent completion of $\AA$ can be realized as the full additive subcategory of the idempotent completion $\widecheck{\AA}$ containing $\AA$ and all kernels of retractions in $\AA$ (see \cite[remark 7.8]{Buhler10} or \cite[proposition~A.11]{HenrardVanRoosmalen19b}).

\begin{remark}\label{remark:ComplementsOfObjectsInWeakIdempotentCompletions}
	Let $X$ be an object of the idempotent completion $\widecheck{\AA}$.  By the construction of the weak idempotent completion $\widehat{\AA}$, we have that $X \in \widehat{\AA}$ if and only if there exists an object $X_c\in \Ob(\AA)$ such that $X\oplus X_c\in \Ob(\AA)$.
\end{remark}

The following observation will be useful.

\begin{lemma}\label{lemma:EmbeddingInWeakIdempotentCompletionExact}
Let $\AA$ be an additive category.  The embedding $i\colon \AA \to \widehat{\AA}$ commutes with limits and colimits.
\end{lemma}

\begin{proof}
Let $D\colon J \to \AA$ be a diagram in $\AA$ and let $L = \lim D.$  For any object $A \in \widehat{A}$, we write $C_A\colon J \to \widehat{A}$ for the constant functor $J \to \AA$ mapping every object to $A$ and every morphism to $1_A$.

For every $X \in \widehat{\AA}$, there is a split kernel-cokernel pair $X \to Y \to Z$ with $Y, Z \in \AA$.  Hence, there is a (split) exact sequence $0 \to \Hom(X,-) \to \Hom(Y,-) \to \Hom(Z,-) \to 0$ in the functor category $\Hom(\widehat{\AA}, \Ab).$  As the natural morphisms $\Hom(Y,L) \to \Hom(C_Y, D)$ and $\Hom(Z,L) \to \Hom(C_Z, D)$ are isomorphisms (as $L$ is the limit of $D$ in $\AA$), we find that the is a natural morphism $\Hom(X,L) \to \Hom(C_X, D)$ is an isomorphism as well.  This shows that $L$ is also the limit of $D$ in $\widehat{\AA}$, as required.

That $i\colon \AA \to \widehat{\AA}$ commutes with colimits is dual.
\end{proof}
\section{Variants of the obscure axiom}

We introduce several variants of the obscure axiom. The variants gain significance throughout this text.

\begin{definition}\label{definition:VariantsOfObscureAxioms}
	Let $\EE$ be a conflation category. We introduce the following axioms.
	\begin{enumerate}[align=left]
		\myitem{\mbox{\textbf{R3}$^-$}}\label{R3-} If the composition $Y'\xrightarrow{f'} Y\xrightarrow{g} Z$ is a deflation and all pullbacks along $g\colon Y\rightarrow Z$ exist, then $g\colon Y\rightarrow Z$ is a deflation.
		\myitem{\mbox{\textbf{R3}$^+$}}\label{R3+} If the composition $Y'\xrightarrow{f'} Y\xrightarrow{g} Z$ is a deflation then $g\colon Y\to Z$ is a deflation.
	\end{enumerate}
	The following axioms are dual.
	\begin{enumerate}[align=left]
		\myitem{\mbox{\textbf{L3}$^-$}}\label{L3-} If a composition $X\xrightarrow{f} Y\xrightarrow{g'} Y'$ is an inflation and all pushouts along $f\colon X\to Y$ exist, then $f\colon X\rightarrow Y$ is an inflation.
		\myitem{\mbox{\textbf{L3}$^+$}}\label{L3+} If a composition $X\xrightarrow{f} Y\xrightarrow{g'} Y'$ is an inflation then $f\colon X\to Y$ is an inflation.
	\end{enumerate}
\end{definition}

\begin{remark}
	For a deflation-exact category, we have
		\[\ref{R3+} \Rightarrow \ref{R3} \Rightarrow \ref{R3-} \Rightarrow \ref{R0*},\]
		meaning that axiom \ref{R3+} is the strongest and axiom \ref{R0*} is the weakest.  A similar sequence exists for inflation-exact categories, changing the axioms to their inflation counterparts.
\end{remark}

The following proposition is a small extension of \cite[proposition~6.4]{BazzoniCrivei13}. It states that, for weakly idempotent complete one-sided exact categories, all aforementioned variants of the obscure axiom are equivalent.

\begin{proposition}\label{proposition:WeaklyIdempotentCompleteR3}
Let $\EE$ be a deflation-exact category. The following are equivalent:
\begin{enumerate}
	\item\label{enumerate:WicR3Plus} $\EE$ satisfies axiom \ref{R3+},
	\item\label{enumerate:WicR3} $\EE$ is weakly idempotent complete and satisfies axiom \ref{R3},
	\item\label{enumerate:WicR3Minus} $\EE$ is weakly idempotent complete and satisfies axiom \ref{R3-}.
\end{enumerate}
\end{proposition}

\begin{proof}
If \eqref{enumerate:WicR3Plus} holds, then \eqref{enumerate:WicR3} holds (and hence also \eqref{enumerate:WicR3Minus}), see \cite[proposition~6.4]{BazzoniCrivei13}.  For the other direction, assume that \eqref{enumerate:WicR3Minus} holds.  Let $g\colon Y \to Z$ be a morphism such that $g\circ f'$ is a deflation (for some $f'\colon Y' \to Y$) as in the statement of axiom \ref{R3+}.  We need to show that $g$ is a deflation.  For this, it suffices to show that $g$ admits all pullbacks.  Equivalently, we can show that, for any $h\colon H \to Z$, the morphism $\begin{psmallmatrix}g & h\end{psmallmatrix}\colon Y \oplus H \to Z$ has a kernel.  We consider the following pullback diagram:
\[\xymatrix{
P \ar@{->>}[d] \ar[rr]^{p} && {Y'} \ar@{->>}[d]^{g \circ f'} \\
{Y \oplus H} \ar[rr]_{\begin{psmallmatrix} g & h \end{psmallmatrix}} && Z}\]
which exists by axiom \ref{R2}.  As $g \circ f' = \begin{psmallmatrix} g & h\end{psmallmatrix} \begin{psmallmatrix}f' \\ 0\end{psmallmatrix}$, the universal property of the pullback shows that $P \to Y'$ is a retraction and thus a deflation by \cref{corollary:RetractionsAreDeflations}. Proposition \ref{proposition:PullbacksPreserveKernels} shows that the composition $\ker p \to P \to Y \oplus H$ is the kernel of $\begin{psmallmatrix} g & h \end{psmallmatrix}\colon Y \oplus H \to Z$.  This completes the proof. 
\end{proof}

The following observation leads to useful equivalent characterizations of the obscure axioms.

\begin{proposition}\label{proposition:WeakFormOfR3}
	Let $\EE$ be a deflation-exact category satisfying axiom \ref{R0*} and let $g\colon Y\to Z$ be a morphism.  The following statements are equivalent:
	\begin{enumerate}
		\item\label{enumerate:EquivalentObscureSetup1} there is a morphism $f'\colon Y'\to Y$ such that the composition $g\circ f'$ is a deflation,
		\item\label{enumerate:EquivalentObscureSetup2} there is an $E \in \EE$ for which $\begin{psmallmatrix}0&g\end{psmallmatrix}\colon E\oplus Y\deflation Z$ is a deflation.
	\end{enumerate}
\end{proposition}

\begin{proof}
Assume that \eqref{enumerate:EquivalentObscureSetup1} holds.  Choose $E=Y'$.	The map $\begin{psmallmatrix}0&g\end{psmallmatrix}$ is obtained as the following composition:
	\[Y'\oplus Y\xrightarrow{\begin{psmallmatrix}1&0\\-f'&1\end{psmallmatrix}}Y'\oplus Y \xrightarrow{\begin{psmallmatrix}gf'&0\\0&1\end{psmallmatrix}}Z\oplus Y\xrightarrow{\begin{psmallmatrix}1&g\\0&1\end{psmallmatrix}}Z\oplus Y\xrightarrow{\begin{psmallmatrix}1&0\end{psmallmatrix}}Z.\] 
	Note that each of the above maps is a deflation by proposition \ref{proposition:BasicProperties} and axiom \ref{R0*}. Axiom \ref{R1} yields the result.
	
	For the converse, it suffices to observe that the deflation $\begin{psmallmatrix}0&g\end{psmallmatrix}\colon E\oplus Y\deflation Z$ factors as $E \oplus Y \to Y \to Z.$
\end{proof}

\begin{proposition}\label{Proposition:ObscureAxiomsCharacterizationsWithSummands}
	Let $\EE$ be a deflation-exact category satisfying axiom \ref{R0*}. 
	\begin{enumerate}
		\item The category $\EE$ satisfies axiom \ref{R3-} if and only if for every morphism $g\colon Y\to Z$ admitting all pullbacks, $\begin{psmallmatrix}0&g\end{psmallmatrix}\colon Y' \oplus Y \to Z$ being a deflation (for any $Y' \in \EE$) implies that $g$ is a deflation.
		\item The category $\EE$ satisfies axiom \ref{R3} if and only if for every morphism $g\colon Y\to Z$ admitting a kernel, $\begin{psmallmatrix}0&g\end{psmallmatrix}\colon Y' \oplus Y \to Z$ being a deflation (for any $Y' \in \EE$) implies that $g$ is a deflation.
		\item The category $\EE$ satisfies axiom \ref{R3+} if and only if for every morphism $g\colon Y\to Z$, if $\begin{psmallmatrix}0&g\end{psmallmatrix}\colon Y' \oplus Y \to Z$ is a deflation (for any $Y' \in \EE$), then $g$ is a deflation.
	\end{enumerate}
\end{proposition}

\begin{proof}
	We show the first equivalence, the other two are similar.
	
	Assume that $\EE$ satisfies axiom \ref{R3-} and let $g\colon Y\to Z$ be a map admitting all pullbacks such that $\begin{psmallmatrix}0&g\end{psmallmatrix}\colon Y'\oplus Y\deflation Z$ is a deflation. As the composition $Y'\oplus Y \xrightarrow{\begin{psmallmatrix}0&1\end{psmallmatrix}} Y\xrightarrow{g}Z$ is a deflation, axiom \ref{R3-} implies that $g$ is a deflation.
	
	Conversely, assume that for every morphism $g\colon Y\to Z$ admitting all pullbacks, $\begin{psmallmatrix}0&g\end{psmallmatrix}$ being a deflation implies that $g$ is a deflation. Let $g\colon Y\to Z$ be a morphism admitting all pullbacks such that there is a map $f'\colon Y'\to Y$ such that $g\circ f'$ is a deflation. By proposition \ref{proposition:WeakFormOfR3}, the map $\begin{psmallmatrix}0&g\end{psmallmatrix}$ is a deflation. By assumption, $g$ is a deflation as well, this shows axiom \ref{R3-}.
\end{proof}

The following proposition is based upon \cite[proposition~6.1]{BazzoniCrivei13} and \cite[lemma~A.2]{ManuelStovicek11}. It states that any deflation-exact category satisfies the dual of theorem \ref{theorem:EquivalentFormulations}.\eqref{item:EquivalentFormulationsC}.

\begin{proposition}\label{proposition:ObscureAxiomOtherSide}
	Let $\EE$ be a deflation-exact category and let $f\colon X\to Y$ be a morphism admitting a cokernel. A map $\begin{psmallmatrix}f\\0\end{psmallmatrix}\colon X \to Y \oplus Y'$ is an inflation if and only if $f$ is an inflation.\\
	If $\EE$ is weakly idempotent complete, one need not require $f$ to have a cokernel.
\end{proposition}

\begin{proof}
	Denote by $g\colon Y\to Z$ the cokernel of $f$ and assume that $\begin{psmallmatrix}f\\0\end{psmallmatrix}\colon X\inflation Y\oplus Y'$ is an inflation. One readily verifies that $\xymatrix{X \ar@{>->}[r]^-{\begin{psmallmatrix}f\\0\end{psmallmatrix}}& Y\oplus Y'\ar@{->>}[r]^-{\begin{psmallmatrix}g&0\\0&1\end{psmallmatrix}} & Z\oplus Y'}$ is a conflation. By axiom \ref{R2} we obtain the following commutative diagram
	\[\xymatrix{
		X\ar@{>->}[r]^-{\begin{psmallmatrix}f\\0\end{psmallmatrix}}\ar@{=}[d] & Y \oplus Y'\ar@{->>}[r]^-{\begin{psmallmatrix}g&0\\0&1\end{psmallmatrix}} & Z\oplus Y'\\
		X\ar@{>->}[r]^f & Y\ar@{->>}[r]^g\ar[u]^{\begin{psmallmatrix}1\\0\end{psmallmatrix}} & Z\ar[u]^{\begin{psmallmatrix}1\\0\end{psmallmatrix}}
	}\] where the lower row is a conflation.  In particular, $f$ is an inflation. 
	
	The converse direction follows from proposition \ref{proposition:BasicProperties}. The last part is \cite[proposition~6.1]{BazzoniCrivei13}.
\end{proof}

\begin{proposition}\label{proposition:R3-DirectSummands}
Let $\EE$ be a deflation-exact category satisfying axiom \ref{R3-}.  Let $f\colon Y \to Z$ and $g\colon Y' \to Z'$ be morphisms such that $\begin{psmallmatrix} f & 0 \\ 0 & g\end{psmallmatrix}\colon Y \oplus Y' \deflation Z \oplus Z'$ is a deflation.  If $f$ admits all pullbacks, then $f$ is a deflation. 
\end{proposition}

\begin{proof}
Axiom \ref{R2} yields a pullback diagram:
\[\xymatrix{
P\ar@{..>>}[d]^{h} \ar@{..>}[r]^{p} & {Y \oplus Y'} \ar@{->>}[d]^{\begin{psmallmatrix} f & 0 \\ 0 & g\end{psmallmatrix}} \\
Z \ar[r]^{\begin{psmallmatrix} 1 \\ 0 \end{psmallmatrix}} & {Z \oplus Z'.}}\]
As ${\begin{psmallmatrix} f & 0 \end{psmallmatrix}} \circ p\colon P \to Z$ is equal to $h$, axiom \ref{R3-} shows that $f$ is a deflation.
\end{proof}
\section{The nine lemma}

Throughout this section, $\EE = (\AA, \bD)$ is a strongly deflation-exact category, i.e.~axiom \ref{R3} holds. The goal of this section is to show that the nine lemma (see theorem \ref{theorem:NineLemma}) holds in $\EE$.  The proof follows \cite{BorceuxBourn04, Bourn01} closely.

\begin{theorem}[The nine lemma]\label{theorem:NineLemma}
	Let $\EE$ be a strongly deflation-exact category. Consider a commutative diagram
	\[\xymatrix{
		X \ar[r]^f\ar@{>->}[d]& Y\ar[r]^g\ar@{>->}[d] & Z\ar@{>->}[d]\\
		X'\ar[r]^{f'}\ar@{->>}[d] & Y'\ar[r]^{g'}\ar@{->>}[d] & Z'\ar@{->>}[d]\\
		X''\ar[r]^{f''} & Y''\ar[r]^{g''} & Z''
	}\] where the columns are conflations and $g'\circ f'=0$. If two of the rows are conflations, so is the third.
\end{theorem}

Before coming to the proof of the nine lemma, we will establish several preliminary results. The first two lemmas are akin to \cite[lemmas~4.2.5 and 4.2.6]{BorceuxBourn04}.

\begin{lemma}\label{lemma:DeflationLemma}
	Consider the following commutative diagram in a strongly deflation-exact category $\EE$
	\[\xymatrix{
		X\ar@{>->}[r]\ar[d]_u & Y\ar@{->>}[r]\ar[d]^v & Z\ar[d]^w\\
		X'\ar@{>->}[r] & Y'\ar@{->>}[r] & Z'\\
	}\] where the rows are conflations.
	\begin{enumerate}
		\item\label{lemma:DeflationLemmaA} If $w$ is an isomorphism, then $u$ is a deflation if and only if $v$ is a deflation.
		\item\label{lemma:DeflationLemmaB} If $u$ and $w$ are deflations then is $v$.
	\end{enumerate}
\end{lemma}

\begin{proof}
	\begin{enumerate}
		\item Assume that $w$ is an isomorphism. If $u$ is a deflation, then \cite[lemma~5.10]{BazzoniCrivei13} shows that $v$ is a deflation as well. 
		
		Conversely, assume that $v$ is a deflation. It follows from \cref{proposition:MitchellPullback} that the left square is a pullback square. Axiom \ref{R2} implies that $u$ is a deflation.
		
		\item Assume that $u$ and $w$ are deflations. Taking the pullback of $Y'\deflation Z'$ along $w$, which exists by axiom \ref{R2}, we obtain the following commutative diagram:
		\[\xymatrix{
			X \ar@{>->}[r]\ar@{->>}[d]_u& Y\ar@{->>}[r]\ar@{.>}[d]_{v'} & Z\ar@{=}[d]\\ 
			X'\ar@{>->}[r]\ar@{=}[d] & P\ar@{->>}[r]\ar@{->>}[d]_{v''} & Z\ar@{->>}[d]_w\\
			X'\ar@{>->}[r] & Y'\ar@{->>}[r] & Z'
		}\] where $v''$ is a deflation and $v'$ is obtained by the universal property of the pullback (thus $v=v''\circ v'$). By the first part of the lemma we find that $v'$ is a deflation. Axiom \ref{R1} now implies that $v$ is a deflation.\qedhere
	\end{enumerate}
\end{proof}

\begin{lemma}\label{lemma:KernelLemma}
	Consider the following commutative diagram in deflation-exact category $\EE$
	\[\xymatrix{
		X \ar@{>->}[r]^f\ar[d]_u& Y\ar@{->>}[r]^g\ar@{->>}[d]_v & Z\ar@{=}[d]\\
		X'\ar@{^{(}->}[r]^{f'} & Y'\ar[r]^{g'} & Z
	}\] where the top row is a conflation and $f'$ is a monomorphism. If $g'\circ f'=0$, then $f'=\ker(g')$.
\end{lemma}

\begin{proof}
	As we have assumed that $f'$ is a monomorphism, we only need to show that $f'$ is a weak kernel of $g'$. To that end, let $t\colon T\to Y'$ be a map such that $g'\circ t=0$. By axiom \ref{R2}, we can take the pullback of $t$ along $v$, obtaining the following commutative solid diagram:
	\[\begin{tikzcd}
	& K\arrow[equal]{r}\arrow[rightarrowtail]{dd}[near start]{k'} & K\arrow[rightarrowtail]{d} & \\
	X\arrow[rightarrowtail, crossing over]{rr}[near start]{f}\arrow[dd, "u"] &&Y\arrow[twoheadrightarrow]{r}{g}\arrow[twoheadrightarrow]{dd}{v} & Z\arrow[equal]{dd}\\
	&P\arrow{ru}{t'}\arrow[twoheadrightarrow]{dd}[near start]{v'}\arrow[lu, "h"', dotted]&&\\
	X'\arrow[hookrightarrow, crossing over]{rr}[near start]{f'}&&Y'\arrow[r, "g'"]&Z\\
	&T\arrow{ru}{t}\arrow[lu, dotted, "l"]&&
\end{tikzcd}\]
Note that $g\circ t'=(g'\circ t)\circ v'=0$ and thus there exists a unique map $h\colon P\to X$ such that $f\circ h=t'$. By commutativity of the diagram, we find that $$t\circ v'=v\circ t'= v\circ f \circ h= f'\circ u\circ h.$$ Since $f'$ is monic and since $f'\circ (u\circ h\circ k') = v \circ t' \circ k' = 0$, we find that $u\circ h\circ k'=0$.  As $v' = \coker (k')$, there exists a unique map $l\colon T\to X'$ such that $l\circ v'=u\circ h$. Thus
\[f' \circ l \circ v' = f'\circ u \circ h = v\circ f \circ h = v \circ t' = t \circ v'.\]
As $v'$ is an epimorphism, we find that $t=f'\circ l$. This shows that $f'$ is a weak  kernel of $g'$.
\end{proof}

The following proposition is a straightforward adaptation of \cite[proposition 1.1.4]{Schneiders99} (see also \cite[corollary 1]{Rump01}).

\begin{proposition}\label{proposition:DeflationMonofactorization}
Let $f\colon X \to Y$ be a morphism in a deflation-exact category $\EE$.  The following are equivalent:
\begin{enumerate}
	\item $f = m \circ e$ where $e$ is a deflation and $m$ is a monomorphism,
	\item $f$ admits a kernel which is an inflation.
\end{enumerate}
In this case, the deflation is given by $e\colon X \deflation X / \ker(f).$
\end{proposition}

\begin{proof}
Assume first that $f = m \circ e$ where $e$ is a deflation and $m$ is a monomorphism.  As $m$ is a monomorphism, we know that $\ker(f) = \ker(e)$.  As $e$ is a deflation, $\ker(f) \to X$ is an inflation.

For the other direction, assume that $i\colon \ker(f) \inflation X$ is an inflation.  Let $C = \coim(f) = X / \ker(f).$  This gives the commutative diagram
\[\xymatrix{
{\ker(f)} \ar@{>->}[r]^-{i} & X \ar[rr]^f \ar@{->>}[dr]_{e} && Y \\ 
&& C \ar[ru]_{m}}\]
We only need to show that $m\colon C \to Y$ is a monomorphism.  To this end, let $t\colon T \to C$ be any morphism such that $m \circ t = 0.$  We need to show that $t=0$.  Using axiom \ref{R2}, we consider the pullback diagram
\[\xymatrix{
P \ar[r]^{t'} \ar@{->>}[d]^{e'} & X \ar@{->>}[d]^{e} \\
T \ar[r]^{t} & C}\]
where the downward arrows are deflations.  As $f \circ t' = m\circ e\circ t' = m\circ t\circ e' = 0$, we know that $t'\colon P \to X$ factors as $t' = i \circ j$ for some $j\colon P \to \ker(f).$  As $e'$ is a deflation (and hence an epimorphism), it now follows from
\[t \circ e' = e \circ t' = e \circ i \circ j = 0\]
that $t = 0.$  This shows that $m$ is a monomorphism.
\end{proof}

We are now in a position the prove the nine lemma.

\begin{proof}[proof of theorem \ref{theorem:NineLemma}]
	\begin{enumerate}
		\item\label{item:MainProofA} It the lower two rows are conflations, so is the third (see \cite[proposition~5.11]{BazzoniCrivei13}).
		\item\label{item:MainProofB} Assume that the upper two rows are conflations. It follows from proposition \ref{proposition:MitchellPullback} that the top-left square is a pullback.  Again using proposition \ref{proposition:MitchellPullback}, we find that $X \inflation X'$ is the kernel of $X' \deflation X'' \xrightarrow{f''} Y''.$  Proposition \ref{proposition:DeflationMonofactorization} yields that $f''\colon X'' \to Y''$ is a monomorphism.
		
		Taking the pullback of $g''$ along the deflation $Z'\deflation Z''$, we obtain the following commutative diagram:
		\[\begin{tikzcd}
			Y\arrow[r, twoheadrightarrow,"g"] \arrow[d, rightarrowtail]& Z\arrow[r, equal]\arrow[d, rightarrowtail] & Z\arrow[d, rightarrowtail]\\
			Y'\arrow[r, dotted] \arrow[d, twoheadrightarrow] &P\arrow[r]\arrow[d, twoheadrightarrow] &Z'\arrow[d, twoheadrightarrow]\\
			Y''\arrow[equal]{r} &Y''\arrow[r, "{g''}"]&Z''
		\end{tikzcd}\]
		Here, the dotted arrow, obtained using the pullback property of $P$, factors the map $g'\colon Y' \to Z'$.  Applying lemma \ref{lemma:DeflationLemma}.\eqref{lemma:DeflationLemmaA} to the first two columns of the previous commutative diagram, we find that the dotted map $Y'\to P$ is a deflation.
		
		Note that $g''\circ f''=0$. Indeed, the composition $X'\xrightarrow{f'}Y'\xrightarrow{g'}Z'\deflation Z''$ equals the composition $X'\deflation X''\xrightarrow{f''}Y''\xrightarrow{g''}Z''$ and is zero as $g'\circ f'=0$, since $X'\deflation X''$ is an epimorphism, we find that $g''\circ f''=0$. 
		
		Using $g''\circ f''=0$ and the pullback property of $P$, we obtain the following commutative diagram:
		\[\begin{tikzcd}
			X'\arrow[r,rightarrowtail,"f'"]\arrow[d, twoheadrightarrow] & Y'\arrow[r, twoheadrightarrow, "g'"]\arrow[d, twoheadrightarrow] & Z'\arrow[d, equal]\\
			X''\arrow[r, dotted]\arrow[d, equal]&P\arrow[r]\arrow[d, twoheadrightarrow]&Z'\arrow[d, twoheadrightarrow]\\
			X''\arrow[r, hookrightarrow, "f''"]&Y''\arrow[r, "g''"]&Z''
		\end{tikzcd}\]
		The dotted map $X''\to P$ is a monomorphism as $f''$ is a monomorphism. By lemma \ref{lemma:KernelLemma}, the map $X''\to P$ is the kernel of the map $P\to Z'$.  As the lower-right square is a pullback, we find that $f''=\ker(g'')$ by proposition \ref{proposition:PullbacksPreserveKernels}. Axiom \ref{R3} now implies that $g''$ is a deflation and thus the lower row is a conflation as well.
		\item\label{item:MainProofC} Assume that the upper and lower rows are conflations. Lemma \ref{lemma:DeflationLemma}.\eqref{lemma:DeflationLemmaB} implies yields that $g'$ is a deflation. Applying part \eqref{item:MainProofA} of this proof to the right columns yields the commutative diagram
		\[\begin{tikzcd}
			X\arrow[r, equal]\arrow[d, rightarrowtail] & X\arrow[r, rightarrowtail, "f"]\arrow[d, rightarrowtail] & Y\arrow[r, twoheadrightarrow, "g"]\arrow[d, rightarrowtail] & Z\arrow[d, rightarrowtail]\\
			X'\arrow[r, dotted]\arrow[d, twoheadrightarrow]&K\arrow[r, rightarrowtail]\arrow[d, twoheadrightarrow] & Y'\arrow[r, twoheadrightarrow, "g'"]\arrow[d, twoheadrightarrow]&Z'\arrow[d, twoheadrightarrow]\\
			X''\arrow[r, equal]&X''\arrow[r, rightarrowtail, "f''"]&Y''\arrow[r, twoheadrightarrow, "g''"]&Z''
		\end{tikzcd}\]
		Here, the dotted arrow is obtained by using that $K=\ker(g')$ and that $g\circ f'=0$.  The short five lemma (the lemma holds for one-sided exact categories, see \cite[lemma~5.3]{BazzoniCrivei13}) implies that $X'\to K$ is an isomorphism.  This concludes the proof. \qedhere
	\end{enumerate}
\end{proof}

\begin{remark}\label{remark:NotHomological}
In \cite{BorceuxBourn04,Bourn01}, it was shown that the nine lemma holds for homological categories.  Although our proof follows these references closely, we do not assume that the category $\EE$ is finitely complete, nor do we assume that every cokernel is a deflation.
\end{remark}
\section{The snake lemma}

In this section we show the snake lemma (see \cref{theorem:SnakeLemmaII}) holds in a deflation-exact category $\EE = (\AA, \bD)$ satisfying axiom \ref{R3+}, i.e.~$\EE$ is weakly idempotent complete and satisfies axiom \ref{R3}.  Our proof follows \cite{Buhler10} closely, and we obtain the snake lemma as a consequence of the ker-coker sequence (see proposition \ref{proposition:Ker-Coker-Sequence} below).

We start with the observation that, in this setting, admissible morphisms are stable under pullbacks along deflations and, when those exist, pushouts along inflations.

\begin{proposition}\label{proposition:AdmissibleStableUnderDiagrams}
	Let $\EE$ be a deflation-exact category satisfying axiom \ref{R3+}. 
		\begin{enumerate}
			\item Given a commutative diagram of the form:
			\begin{gather}\begin{aligned}\xymatrix{
				A_1\ar@{>->}[r]\ar@{=}[d] & A_2\ar[d]^{f_2}\ar@{->>}[r] & A_3\ar[d]^{f_3}\\
				A_1\ar@{>->}[r] & B_2\ar@{->>}[r] & B_3
			}\end{aligned}\label{DiagramEqualLeft}\end{gather}
			The map $f_2$ is admissible if and only if $f_3$ is admissible.
				\item Given a commutative diagram of the form:
				\begin{gather}\begin{aligned}\xymatrix{
				A_1\ar@{>->}[r]\ar[d]^{f_1} & A_2\ar[d]^{f_2}\ar@{->>}[r] & A_3\ar@{=}[d]\\
				B_1\ar@{>->}[r] & B_2\ar@{->>}[r] & A_3
			}\end{aligned}\label{DiagramEqualRight}\end{gather}
			The map $f_1$ is admissible if and only if $f_2$ is admissible.
		\end{enumerate}
\end{proposition}

\begin{proof}
	\begin{enumerate}
		\item Note that the right square of diagram \ref{DiagramEqualLeft} is both a pullback square and a pushout square. Indeed this follows from \cite[proposition~3.7]{HenrardVanRoosmalen19a}).
		
		Now assume that $f_2$ is admissible.  It follows from proposition \ref{proposition:PullbacksPreserveKernels} that $\coker(f_2) \cong \coker(f_3).$  Axiom \ref{R3+} implies that the induced map $B_3\to \coker(f_2)$ is a deflation. Hence we obtain the following commutative diagram:
				\[\xymatrix{
					A_1\ar@{.>}[r]\ar@{=}[d] & \im(f_2)\ar@{>->}[d]\ar@{.>}[r] & K\ar@{>->}[d]\\
					A_1\ar@{>->}[r]\ar@{->>}[d] & B_2\ar@{->>}[r]\ar@{->>}[d] & B_3\ar@{->>}[d]\\
					0\ar@{>->}[r] & \coker(f_2)\ar@{=}[r] & \coker(f_2)
				}\]
				By the nine lemma, the upper row is a conflation as well. Using that $K$ is a kernel, we find that $f_3$ factors through $K\inflation B_3$. The induced map $A_3\to K$ is a deflation by axiom \ref{R3+} and by commutativity of the square $A_2A_3K\im(f_2)$. Indeed, one can verify the commutativity of the square $A_2A_3K\im(f_2)$ by using that $K\inflation B_3$ is a monomorphism and that the square $A_2A_3B_3B_2$ commutes. This shows that $f_3$ is admissible and that $K\cong \im(f_3)$.
				
		Conversely, assume that $f_3$ is admissible.  It follows from \cite[proposition~3.10]{HenrardVanRoosmalen19a} and axiom \ref{R3} that admissible maps are stable under pullbacks along deflations. Note that this does not require $\EE$ to be weakly idempotent complete.
		
		\item The left square of diagram \ref{DiagramEqualRight} is both a pullback square and a pushout square by \cite[proposition~3.7]{HenrardVanRoosmalen19a}).
		
		Assume that $f_1$ admissible. As the left square is a pushout, \cref{proposition:PullbacksPreserveKernels} shows that $\coker(f_1) \cong \coker(f_2).$  Axiom \ref{R3+} implies that the map $B_2\to \coker(f_1)$ is a deflation.  Hence, we obtain the following commutative diagram:
				\[\xymatrix{
					\im(f_1)\ar@{.>}[r]\ar@{>->}[d] & K\ar@{.>}[r]\ar@{>->}[d] & A_3\ar@{=}[d]\\
					B_1\ar@{->>}[r]\ar@{->>}[d] & B_2\ar@{->>}[r]\ar@{->>}[d] & A_3\ar@{->>}[d]\\
					\coker(f_1)\ar@{=}[r] & \coker(f_1)\ar@{->>}[r] & 0
				}\]
				By the nine lemma, the upper row is a conflation. Note that $f_2$ factors through $K\inflation B_2$, moreover, we obtain the following commutative diagram:
				\[\xymatrix{
					\ker(f_1)\ar@{=}[r]\ar@{>->}[d] & \ker(f_1)\ar@{->>}[r]\ar[d] & 0\ar@{>->}[d]\\
					A_1\ar@{>->}[r]\ar@{->>}[d] & A_2\ar@{->>}[r]\ar@{.>}[d] & A_3\ar@{=}[d]\\
					\im(f_1)\ar@{>->}[r] & K\ar@{->>}[r] & A_3
				}\]
				Applying the nine lemma once more, we find that the middle column is a conflation. It follows that $f_2$ is admissible and that $K\cong \im(f_2)$.
		
		Conversely, assume that $f_2$ is admissible. It is straightforward to see that $A_2\deflation A_3$ factors through $A_2\deflation \im(f_2)$. By axiom \ref{R3+}, the induced map $\im(f_2)\to A_3$ is a deflation. Using that the left square is a pullback, we obtain the following commutative diagram:
				\[\xymatrix{
					\ker(f_1)\ar@{=}[r]\ar[d] & \ker(f_1)\ar@{->>}[r]\ar@{>->}[d] & 0\ar@{>->}[d]\\
					A_1\ar@{>->}[r]\ar[d] & A_2\ar@{->>}[r]\ar@{->>}[d] & A_3\ar@{=}[d]\\
					K\ar@{>->}[r] & \im(f_2)\ar@{->>}[r] & A_3
				}\]
				By the nine lemma, the left column is a conflation. Clearly, $f_1$ factors through $A_1\deflation K$. It remains to show that the induced map $K\to B_1$ is an inflation. Applying the nine lemma once more to the commutative diagram
				\[\xymatrix{
					K\ar[d]\ar@{>->}[r] & \im(f_2)\ar@{->>}[r]\ar@{>->}[d] & A_3\ar@{=}[d]\\
					B_1\ar@{>->}[r]\ar[d] & B_2\ar@{->>}[r]\ar@{->>}[d] & A_3\ar@{->>}[d]\\
					\coker(f_2)\ar@{=}[r] & \coker(f_2)\ar@{->>}[r] & 0
				}\] we find that the left column is a conflation. It follows that $f_1$ is admissible and that $K\cong \im(f_1)$. \qedhere
	\end{enumerate}
\end{proof}

The next proposition shows that the ker-coker-sequence property holds for weakly idempotent complete deflation-exact categories. We generalize the proof given in \cite[proposition~8.11]{Buhler10}.

\begin{proposition}[Ker-coker-sequence]\label{proposition:Ker-Coker-Sequence}
	Let $\EE$ be a deflation-exact category satisfying axiom \ref{R3+}. Let $\xymatrix{A\ar[r]|{\circ}^f&B }$ and $\xymatrix{B\ar[r]|{\circ}^g&C}$ be admissible morphisms such that $h=g\circ f$ is admissible as well. There is a natural exact sequence 
	\[0 \to \ker(f) \to \ker(h) \to \ker(g) \to \coker(f) \to \coker(h) \to \coker(g)\to 0.\]
\end{proposition}

\begin{proof}
	Clearly $g$ induces a map $\im(f)\to \im(h)$, moreover, by axiom \ref{R3+}, this map is a deflation. Hence, we obtain a conflation $X\inflation \im(f)\deflation \im(h)$. Consider the following commutative diagram:
	\[\xymatrix{
		\ker(f) \ar@{.>}[r]\ar@{=}[d]& \ker(h)\ar@{.>}[r]\ar@{>->}[d] & X\ar@{>->}[d]\\
		\ker(f)\ar@{>->}[r]\ar@{->>}[d] & A\ar@{->>}[r]\ar@{->>}[d] & \im(f)\ar@{->>}[d]\\
		0\ar@{>->}[r] & \im(h)\ar@{=}[r] & \im(h)
	}\]
	By the nine lemma, the upper row of this diagram is a conflation.
	
	Similarly, $f$ induces a map $\coker(h)\to \coker(g)$ which is a deflation by axiom \ref{R3+}. Hence, we obtain a conflation $Z\inflation \coker(h)\deflation \coker(g)$. Consider the following commutative diagram:
	\[\xymatrix{
		\im(h)\ar@{.>}[r]\ar@{=}[d] & \im(g)\ar@{.>}[r]\ar@{>->}[d] & Z\ar@{>->}[d]\\
		\im(h)\ar@{>->}[r]\ar@{->>}[d] & C\ar@{->>}[r]\ar@{->>}[d] & \coker(h)\ar@{->>}[d]\\
		0\ar@{>->}[r] & \coker(g) \ar@{=}[r] & \coker(g)
	}\]
	Again, the nine lemma yields that the top row is a conflation.
	
	By axiom \ref{R1}, the composition $B\deflation \im(g)\deflation Z$ is a deflation. Axiom \ref{R3+} implies that the induced map $\coker(f)\to Z$ such that the square $B\coker(f)Z\im(g)$ is commutative, is a deflation. Hence we obtain a conflation $Y\inflation \coker(f)\deflation Z$. Consider the following natural commutative diagram:
	\[\xymatrix{
		X\ar@{.>}[r]\ar@{>->}[d] & \ker(g)\ar@{.>}[r]\ar@{>->}[d] & Y\ar@{>->}[d]\\
		\im(f)\ar@{>->}[r]\ar@{->>}[d] & B\ar@{->>}[r]\ar@{->>}[d] & \coker(f)\ar@{->>}[d]\\
		\im(h)\ar@{>->}[r] & \im(g)\ar@{->>}[r] & Z
	}\]
	Again, the nine lemma implies that the top row of this diagram is a conflation. The result follows by gluing the edges of the previous three diagrams together.
\end{proof}

As a corollary, we find obtain the short snake lemma.

\begin{corollary}[Short snake lemma]\label{corollary:SnakeLemmaI}
	Let $\EE$ be a deflation-exact category satisfying axiom \ref{R3+}. Consider a commutative diagram 
	\[\xymatrix{
		A_1\ar@{>->}[r]^{\phi_1}\ar[d]|{\circ}^{f_1} & A_2\ar@{->>}[r]^{\phi_2}\ar[d]|{\circ}^{f_2} & A_3\ar[d]|{\circ}^{f_3}\\
		B_1\ar@{>->}[r]_{\phi_1'} & B_2\ar@{->>}[r]_{\phi_2'} & B_3
	}\]
	such that the rows are conflations and the $f_i$'s are admissible morphisms. Write $K_i=\ker(f_i)$ and $C_i=\coker(f_i)$. There exists a natural connecting morphism $\delta$ such that the sequence 
	\[\xymatrix{
		0\ar[r]&K_1\ar@{>->}[r]^{\psi_1} & K_2\ar[r]^{\psi_2} & K_3\ar[r]^{\delta} & C_1\ar[r]^{\psi_1'} & C_2\ar@{->>}[r]^{\psi_2'} & C_3\ar[r] & 0
	}\] is an exact sequence. Here the maps $\psi_i$ and $\psi_i'$ are the natural maps induced by the above diagram.
\end{corollary}

\begin{proof}
	By \cite[proposition~5.2]{BazzoniCrivei13}, the given map between conflations factors as 
	\[\xymatrix{
		A_1\ar@{>->}[r]^{\phi_1}\ar[d]|{\circ}^{f_1} & A_2\ar@{->>}[r]^{\phi_2}\ar[d]|{\circ}^{f_2'} & A_3\ar@{=}[d]\\
		B_1\ar@{=}[d]\ar@{>->}[r] & P\ar@{->>}[r]\ar[d]|{\circ}^{f_2''} & A_3\ar[d]|{\circ}^{f_3}\\
		B_1\ar@{>->}[r]_{\phi_1'} & B_2\ar@{->>}[r]_{\phi_2'} & B_3	
	}\] such that $f_2=f_2''\circ f_2'$, the upper-left and lower-right squares are bicartesian squares, i.e.~both pullbacks and pushouts. By proposition \ref{proposition:AdmissibleStableUnderDiagrams}, both $f'_2$ and $f''_2$ are admissible.
	
	Applying proposition \ref{proposition:Ker-Coker-Sequence} to the composition $f_2=f_2''\circ f_2'$, yields the exact sequence
	\[0 \to \ker(f_2') \to \ker(f_2) \to \ker(f_2'') \to \coker(f_2') \to \coker(f_2) \to \coker(f_2'')\to 0.\]
	Using that pullbacks preserve kernels, pushouts preserve cokernels (proposition \ref{proposition:PullbacksPreserveKernels}), yields the required sequence.
\end{proof}

\begin{remark}
	\begin{enumerate}
		\item The short snake lemma holds for a strongly one-sided exact category $\EE$ only if $\EE$ is weakly idempotent complete. The proof is similar to \cite[remark~8.14]{Buhler10}.
		\item The proof of proposition \ref{proposition:Ker-Coker-Sequence} is similar to the proof given in \cite[proposition~8.11]{Buhler10}. However, we avoid using axiom \ref{L3} by combining axiom \ref{R3} and the nine lemma to obtain the desired inflations.
	\end{enumerate}
\end{remark}

We now come to the full version of the snake lemma. 

\begin{theorem}[Snake lemma]\label{theorem:SnakeLemmaII}
	Let $\EE$ be a deflation-exact category satisfying axiom \ref{R3+}. Consider a commutative diagram
	\begin{gather}\begin{aligned}\xymatrix{
		& A_1\ar[r]|{\circ}^{\phi_1}\ar[d]|{\circ}^{f_1} & A_2\ar[r]|{\circ}^{\phi_2}\ar[d]|{\circ}^{f_2} & A_3\ar@{->>}[r]\ar[d]|{\circ}^{f_3} &0 \\
		0\ar@{>->}[r] & B_1\ar[r]|{\circ}_{\phi_1'} & B_2\ar[r]|{\circ}_{\phi_2'} & B_3
	}\end{aligned}\label{diagram:SnakeDiagramII}\end{gather}
	with exact rows. There is an induced exact sequence
	\[\xymatrix{\ker(\phi_1)\ar@{>->}[r] & K_1\ar[r] & K_2\ar[r] & K_3\ar[r] & C_1\ar[r] & C_2\ar[r] & C_3\ar@{->>}[r] & \coker(\phi_2')}\]
	where $K_i = \ker(f_i)$ and $C_i = \coker(f_i)$.
\end{theorem}

\begin{proof}
	Note that $\phi_2$ is a deflation, $\phi_1'$ is an inflation, and that the composition $\ker(\phi_1)\inflation A_1\deflation \im(f_1)$ is zero. By axiom \ref{R3+}, the induced map $\im(\phi_1)\to \im(f_1)$ is a deflation. By the nine lemma, the top row of the following commutative diagram is a conflation:
	\[\xymatrix{
		\ker(\phi_1)\ar@{.>}[r]\ar@{=}[d] & \ker(f_1)\ar@{.>}[r]\ar@{>->}[d] & K\ar@{>->}[d]\\
		\ker(\phi_1)\ar@{>->}[r]\ar@{->>}[d] & A_1\ar@{->>}[r]\ar@{->>}[d] & \im(\phi_1)\ar@{->>}[d]\\
		0\ar@{>->}[r] & \im(f_1)\ar@{=}[r] & \im(f_1)
	}\]
	
	By axiom \ref{R1}, the composition $A_2\deflation A_3\deflation \im(f_3)$ is a deflation and thus an epimorphism. It follows that the composition $\im(f_3)\inflation B_3 \deflation \coker(\phi_2')$ is zero if and only if the natural map $A_2\to \coker(\phi_2')$ is zero. The latter follows from the commutativity of diagram \eqref{diagram:SnakeDiagramII}. Using axiom \ref{R3+}, we find that the induced map $\coker(f_3)\to \coker(\phi_2')$ is a deflation. The nine lemma yields that the left column of the following commutative diagram is a conflation:
	\[\xymatrix{
		\im(f_3)\ar@{=}[r]\ar@{.>}[d] & \im(f_3)\ar@{->>}[r]\ar@{>->}[d] & 0\ar@{>->}[d]\\
		\im(\phi_2')\ar@{>->}[r]\ar@{.>}[d] & B_3\ar@{->>}[r]\ar@{->>}[d] & \coker(\phi_2')\ar@{=}[d]\\
		L\ar@{>->}[r] & \coker(f_3)\ar@{->>}[r] & \coker(\phi_2')
	}\]
	
	Applying corollary \ref{corollary:SnakeLemmaI} to the commutative diagram
	\[\xymatrix{
		\im(\phi_1)\ar@{>->}[rr]\ar@{->>}[d] && A_2\ar@{->>}[rr]\ar@{->>}[d] && A_3\ar@{->>}[d]\\
		\im(f_1)\ar@{>->}[d] &&\im(f_2)\ar@{>->}[d] && \im(f_3)\ar@{>->}[d]\\
		B_1\ar@{>->}[rr] && B_2\ar@{->>}[rr] && \im(\phi_2')
	}\] together with the above yields the desired result.
\end{proof}
\section{Equivalent formulations of the obscure axiom}\label{section:EquivalentFormulations}

In this section, we expand upon the previous results and provide several equivalent formulations of the obscure axiom for a deflation- or inflation-exact category $\EE$ (theorems \ref{theorem:EquivalentFormulations} and \ref{theorem:WeaklyIdempotentCompleteEquivalentFormulations}). In particular, we show some converses of the previous sections, i.e.~we show that the obscure axiom holds in $\EE$ if and only if the nine lemma holds, and that axiom \ref{R3+} holds if and only if the snake lemma or the ker-coker sequence property holds.

\subsection{Equivalent formulations of axiom \texorpdfstring{\ref{R3}}{R3}}

We now come to the proof of theorem \ref{theorem:EquivalentFormulations}.  For clarity, the proof has been split in several parts.

\begin{proof}[Proof of equivalences $\eqref{item:EquivalentFormulationsA} \Leftrightarrow \eqref{item:EquivalentFormulationsB}$] 	Theorem \ref{theorem:NineLemma} shows the implication $\eqref{item:EquivalentFormulationsA}\Rightarrow\eqref{item:EquivalentFormulationsB}$.  For the implication $\eqref{item:EquivalentFormulationsB}\Rightarrow\eqref{item:EquivalentFormulationsA}$, consider the commutative diagram
	\begin{equation}\label{equation:DiagramR3}
		\begin{tikzcd}
			&&L\arrow[d, rightarrowtail]\\
			&&A\arrow[d, twoheadrightarrow]\arrow[dl, "i"']\\
			K\arrow[r, "k"]&B\arrow[r, "p"]&C
		\end{tikzcd}
	\end{equation} 
	where $k=\ker(p)$. We need to show that $p$ is a deflation. By axiom \ref{R2}, we can consider the pullback $P$ of $p$ along the deflation $p\circ i$.  It follows from the universal property of the pullback that the morphism $P \to A$ is a retraction, moreover, by \cref{proposition:PullbacksPreserveKernels} it has a kernel and thus \cref{proposition:AxiomR0*} yields that it is a deflation and that $P\cong K\oplus A$. We obtain the following commutative diagram:
	\begin{gather}\begin{aligned}
		\xymatrix{
			0 \ar@{>->}[r]\ar@{>->}[d]& L\ar@{=}[r]\ar@{>->}[d] & L\ar@{>->}[d]\\
			K\ar@{>->}[r]\ar@{=}[d] & K\oplus A\ar@{->>}[r]\ar@{->>}[d]_{\begin{psmallmatrix}k& i\end{psmallmatrix}} & A\ar@{->>}[d]\ar@{.>}[dl]^i\\
			K\ar[r]_k & B\ar[r]_p & C
		}
	\end{aligned}\label{equation:DiagramR3Pullback}\end{gather}
	By the nine lemma, the lower row is a conflation and thus $p$ is a deflation. We conclude that $\eqref{item:EquivalentFormulationsA}\Leftrightarrow \eqref{item:EquivalentFormulationsB},$ as required. \end{proof}
	
\begin{proof}[Proof of equivalences $\eqref{item:EquivalentFormulationsA} \Leftrightarrow \eqref{item:EquivalentFormulationsE} \Leftrightarrow \eqref{item:EquivalentFormulationsF}\Leftrightarrow \eqref{item:EquivalentFormulationsC}$]

The equivalence $\eqref{item:EquivalentFormulationsA}\Leftrightarrow\eqref{item:EquivalentFormulationsC}$ is shown in proposition \ref{Proposition:ObscureAxiomsCharacterizationsWithSummands}. The implication $\eqref{item:EquivalentFormulationsA}\Rightarrow \eqref{item:EquivalentFormulationsE}$ is shown in \cite[proposition~2.6]{HenrardVanRoosmalen19b}. The implication $\eqref{item:EquivalentFormulationsE}\Rightarrow \eqref{item:EquivalentFormulationsF}$ is trivial.  To show the implication $\eqref{item:EquivalentFormulationsF}\Rightarrow \eqref{item:EquivalentFormulationsA}$, we consider the setup given in the commutative diagram \eqref{equation:DiagramR3}, with $k=\ker(p)$.  As in the proof following diagram \eqref{equation:DiagramR3}, we consider the pullback $P$ of $p$ along the deflation $p\circ i$ and obtain the lower-right square of \eqref{equation:DiagramR3Pullback} as the pullback diagram.  Using \cite[proposition~5.7]{BazzoniCrivei13}, we find the conflation
\[\xymatrix{
{K \oplus A} \ar@{>->}[rr]^{\begin{psmallmatrix} k & i \\ 0 & -1 \end{psmallmatrix} } &&
{B \oplus A} \ar@{->>}[rr]^{\begin{psmallmatrix} p & pi \end{psmallmatrix}} && C }\]
We can now find $K \to B \to C$ as a direct summand of this conflation.  Hence, it follows from \eqref{item:EquivalentFormulationsF} that $K \to B \to C$ is a conflation.  This shows the implication $\eqref{item:EquivalentFormulationsF}\Rightarrow \eqref{item:EquivalentFormulationsA}$.
\end{proof}

\begin{proof}[Proof of equivalences $\eqref{item:EquivalentFormulationsA}\Leftrightarrow \eqref{item:EquivalentFormulationsD}$]
	The implication $\eqref{item:EquivalentFormulationsA}\Rightarrow \eqref{item:EquivalentFormulationsD}$ is trivial. Conversely, assume $\eqref{item:EquivalentFormulationsD}$ and let $g\colon Y\to Z$, $f\colon X\to Y$ and $f'\colon Y'\to Y$ be maps such that $g\circ f'$ is a deflation and $f=\ker(g)$. By axiom \ref{R2}, we obtain the commutative diagram
	\[\xymatrix{
		X\ar[r]^{i}\ar@{=}[d] & P\ar[r]^{p}\ar@{->>}[d]^h & Y'\ar@{->>}[d]^{gf'}\\
		X\ar[r]^f & Y\ar[r]^g & Z
	}\]
	where the right square is a pullback. By the pullback property, there is a unique map $u\colon Y'\to P$ such that $gf'=hu$ and $pu=1_{Y'}$. It follows that $p$ is a retraction with kernel $i$. By proposition \ref{proposition:AxiomR0*} and axiom \ref{R1}, the composition $gh=gf'p$ is a deflation. By $\eqref{item:EquivalentFormulationsD}$, $g$ is a deflation as required.
\end{proof}

\begin{proof}[Proof of equivalences $\eqref{item:EquivalentFormulationsA} \Leftrightarrow \eqref{item:EquivalentFormulationsH} \Leftrightarrow \eqref{item:EquivalentFormulationsI}$]

The equivalence of $\eqref{item:EquivalentFormulationsA}\Leftrightarrow \eqref{item:EquivalentFormulationsH}$ is shown in \cite[proposition~6.2]{HenrardVanRoosmalen19b}.  The equivalence of $\eqref{item:EquivalentFormulationsH}\Leftrightarrow \eqref{item:EquivalentFormulationsI}$ follows from \cite[theorem~1.2 and proposition~6.2]{HenrardVanRoosmalen19b}.
\end{proof}

\begin{proof}[Proof of equivalences $\eqref{item:EquivalentFormulationsA} \Leftrightarrow \eqref{item:EquivalentFormulationsJ}$]

	Assume that \eqref{item:EquivalentFormulationsA} holds, we need to show that $\EE$ lies extension closed in its weak idempotent completion. Let $X\stackrel{i}{\inflation} Y\stackrel{p}{\deflation} Z$ be a conflation in $\widehat{\EE}$ such that $X,Z\in \EE$. By definition, there exists a conflation $X'\stackrel{i'}{\inflation} Y'\stackrel{p'}{\deflation} Z'$ in $\widehat{\EE}$ such that the direct sum $X\oplus X'\inflation Y\oplus Y'\deflation Z\oplus Z'$ is a conflation in $\EE$. Without loss of generality, we may assume that $X'\inflation Y'\deflation Z'$ is a conflation in $\EE$. By remark \ref{remark:ComplementsOfObjectsInWeakIdempotentCompletions}, there exists an object $Y_c\in \EE$ such that $Y\oplus Y_c\in \EE$. It follows that the following natural commutative diagram lies in $\EE$:
	\[\xymatrix{
		X\oplus X'\ar@{>->}[r]^-{\begin{psmallmatrix}i&0\\0&i'\\0&0\end{psmallmatrix}} & Y\oplus Y'\ar@{->>}[r]^-{\begin{psmallmatrix}p&0&0\\0&p'&0\\0&0&1\end{psmallmatrix}}\ar[d]^-{\begin{psmallmatrix}1&0&0\\0&p'&0\\0&0&1\end{psmallmatrix}} \oplus Y_c & Z\oplus Z'\oplus Y_c\ar@{=}[d]\\
		X \ar[r]_-{\begin{psmallmatrix}i\\0\\0\end{psmallmatrix}}& Y\oplus Z'\oplus Y_c\ar[r]_-{\begin{psmallmatrix}p&0&0\\0&1&0\\0&0&1\end{psmallmatrix}} & Z\oplus Z'\oplus Y_c
	}\]
	Note that the lower row of the above diagram is a kernel-cokernel pair in $\EE$. By axiom \ref{R3}, the lower row is a conflation in $\EE$. By axiom \ref{R2}, pullbacks along $Y\oplus Z'\oplus Y_c\to Z\oplus Z'\oplus Y_C$ exist in $\EE$:
	\[\xymatrix{
Y\oplus Z'\oplus Y_c\ar@{->>}[r] & Z\oplus Z'\oplus Y_c\\
P\ar@{->>}[r]\ar[u]^-{\begin{psmallmatrix}1\\0\\0\end{psmallmatrix}} & Z\ar[u]^-{\begin{psmallmatrix}1\\0\\0\end{psmallmatrix}}
	}\]
	It now follows from lemma \ref{lemma:EmbeddingInWeakIdempotentCompletionExact} that $P \cong Y$.  Hence, $Y \in \EE$.  This shows that $\EE$ is extension-closed in $\widehat{\EE}.$
	
	Conversely, assume that \eqref{item:EquivalentFormulationsJ} holds.  Let $X\xrightarrow{f}Y\xrightarrow{g} Z$ be a kernel-cokernel pair in $\EE$ and let $h\colon A\to Y$ be map such that $g\circ h$ is a deflation in $\EE$. Since $g\circ h$ is also a deflation in $\widehat{\EE}$ and $\widehat{\EE}$ satisfies axiom \ref{R3+}, the map $g$ is a deflation in $\widehat{\EE}$. Since $\EE\hookrightarrow \widehat{\EE}$ is fully exact, $Y\xrightarrow{g} Z$ is a deflation in $\EE$ as well.
	\end{proof}

\begin{remark}
The Gabriel-Quillen embedding of an exact category in an abelian category allows one to reduce homological statements in an exact category to homological statements in an abelian category.  The equivalence $\eqref{item:EquivalentFormulationsA}\Leftrightarrow \eqref{item:EquivalentFormulationsI}$ in theorem \ref{theorem:EquivalentFormulations} gives a similar reduction: here, one can establish statements about conflations in a strongly deflation-exact category by reducing it to a similar statement in its exact hull.  In particular, using the equivalence $\eqref{item:EquivalentFormulationsA} \Leftrightarrow \eqref{item:EquivalentFormulationsI}$ in theorem \ref{theorem:EquivalentFormulations}, one can deduce the nine lemma in a strong deflation-exact category $\EE$ from the nine lemma in the exact hull $\overline{\EE}$ of $\EE$.
\end{remark}

\subsection{\texorpdfstring{Equivalent formulations of axiom \ref{R3+}}{Equivalent formulations of axiom R3+}}

We now come to the proof of theorem \ref{theorem:WeaklyIdempotentCompleteEquivalentFormulations}.  We will again present the proof in several parts.

\begin{proof}[Proof of equivalences $\eqref{item:WeaklyIdempotentCompleteEquivalentFormulationsA} \Leftrightarrow \eqref{item:WeaklyIdempotentCompleteEquivalentFormulationsB}\Leftrightarrow \eqref{item:WeaklyIdempotentCompleteEquivalentFormulationsC}$]

This is shown in proposition \ref{proposition:WeaklyIdempotentCompleteR3} and proposition \ref{Proposition:ObscureAxiomsCharacterizationsWithSummands}.
\end{proof}

\begin{proof}[Proof of equivalences $\eqref{item:WeaklyIdempotentCompleteEquivalentFormulationsA} \Leftrightarrow \eqref{item:WeaklyIdempotentCompleteEquivalentFormulationsF}$]

The implication $\eqref{item:WeaklyIdempotentCompleteEquivalentFormulationsA}\Rightarrow\eqref{item:WeaklyIdempotentCompleteEquivalentFormulationsF}$ is trivial.  For the implication $\eqref{item:WeaklyIdempotentCompleteEquivalentFormulationsF}\Rightarrow\eqref{item:WeaklyIdempotentCompleteEquivalentFormulationsA},$ consider morphisms $f\colon X \to Y$ and $g\colon Y \to Z$ such that $g \circ f$ is a deflation.  We find the following pullback square
	\[\xymatrix{
	P \ar[r]^p \ar@{->>}[d]_{h} & X \ar@{->>}[d]^{gf} \ar[ld]^{f} \\
	Y \ar[r]_{g} & Z}\]
		where $h\colon P \deflation Y$ is a deflation by axiom \ref{R2}. Note that $p\colon P\to X$ is a retraction by the pullback property. Write $s\colon X\to P$ for the corresponding section. By axiom \ref{R0*}, the maps $\begin{psmallmatrix}p&1\end{psmallmatrix}$ and $\begin{psmallmatrix}1&s\end{psmallmatrix}$ are deflations. The following commutative diagram
		\[\xymatrix{
			P\oplus X\ar@{->>}[r]^{\begin{psmallmatrix}1&s\end{psmallmatrix}}\ar@{->>}[rd]_{\begin{psmallmatrix}p&1\end{psmallmatrix}} & P\ar[d]^p\\
			& X
		}\] and \eqref{item:WeaklyIdempotentCompleteEquivalentFormulationsF} yields that $p\colon P\to X$ is a deflation. It now follows from axiom \ref{R1} that the composition $g \circ h$ is a deflation.  Hence \eqref{item:WeaklyIdempotentCompleteEquivalentFormulationsF} implies that $g$ is a deflation.  This shows that \eqref{item:WeaklyIdempotentCompleteEquivalentFormulationsA} holds, as required.
\end{proof}

\begin{proof}[Proof of equivalences $\eqref{item:WeaklyIdempotentCompleteEquivalentFormulationsD} \Leftrightarrow \eqref{item:WeaklyIdempotentCompleteEquivalentFormulationsA}  \Leftrightarrow \eqref{item:WeaklyIdempotentCompleteEquivalentFormulationsE}$]

		The implications $\eqref{item:WeaklyIdempotentCompleteEquivalentFormulationsA}\Rightarrow\eqref{item:WeaklyIdempotentCompleteEquivalentFormulationsD}$ and $\eqref{item:WeaklyIdempotentCompleteEquivalentFormulationsA}\Rightarrow\eqref{item:WeaklyIdempotentCompleteEquivalentFormulationsE}$ follow from proposition \ref{proposition:Ker-Coker-Sequence}, corollary \ref{corollary:SnakeLemmaI} and theorem \ref{theorem:SnakeLemmaII}.
		
		We now show that $\eqref{item:WeaklyIdempotentCompleteEquivalentFormulationsD} \Rightarrow \eqref{item:WeaklyIdempotentCompleteEquivalentFormulationsF}$.  As $\eqref{item:WeaklyIdempotentCompleteEquivalentFormulationsF} \Rightarrow \eqref{item:WeaklyIdempotentCompleteEquivalentFormulationsA}$, this is sufficient.  Let $g\colon Y \to Z$ be a morphism and $f\colon X \deflation Y$ be a deflation such that $g \circ f\colon X \deflation Z$ is a deflation, as in the statement of \eqref{item:WeaklyIdempotentCompleteEquivalentFormulationsF}.  We need to show that $g$ is a deflation.  Applying the ker-coker sequence to the composition $\ker f \inflation X \stackrel{gf}{\deflation} Z$ gives an acyclic sequence
		\[\xymatrix@1{0 \ar[r] & X \ar@{>->}[r] & {\ker(gf)} \ar[r]|-{\circ} & Y \ar@{->>}[r] & Z \ar[r] & 0.}\]
	In particular, $g\colon Y \deflation Z$ is a deflation.  The implication $\eqref{item:WeaklyIdempotentCompleteEquivalentFormulationsE} \Rightarrow \eqref{item:WeaklyIdempotentCompleteEquivalentFormulationsF}$ can be shown in a similar fashion.
\end{proof}

\begin{proof}[Proof of implications $\eqref{item:WeaklyIdempotentCompleteEquivalentFormulationsF} \Rightarrow \eqref{item:WeaklyIdempotentCompleteEquivalentFormulationsH}\Rightarrow \eqref{item:WeaklyIdempotentCompleteEquivalentFormulationsG} \Rightarrow \eqref{item:WeaklyIdempotentCompleteEquivalentFormulationsF}$]

Only the last implication is nontrivial, so let $f\colon X \to Y$ be any map and let $g\colon Y \to Z$ be a deflation such that $g \circ f$ is a deflation.  It follows from \cite{HenrardVanRoosmalen19a} (or the dual of \cite[proposition 5.7]{BazzoniCrivei13}) that $\begin{psmallmatrix} gf & f \end{psmallmatrix}\colon X \oplus Y \deflation Z$ is a deflation (this uses axiom \ref{R0*}).  As invertible maps are deflations, axiom \ref{R2} gives a deflation $\begin{psmallmatrix}gf & f\end{psmallmatrix} \begin{psmallmatrix} 1 & 0 \\ -h & 1 \end{psmallmatrix} = \begin{psmallmatrix} 0 & f \end{psmallmatrix} \colon X \oplus Y \deflation Z.$  It now follows from \eqref{item:WeaklyIdempotentCompleteEquivalentFormulationsG} that $f$ is a deflation, as required.
\end{proof}
\section{Obscure closures of one-sided exact categories}\label{section:ObscureClosuresOfOneSidedExactCategories}

Let $\EE = (\AA, \bD)$ be a deflation-exact category satisfying axiom \ref{R0*}.  In this section, we examine closures of $\EE$ under the obscure axiom \ref{R3} (as well as the closures under axioms \ref{R3-} and \ref{R3+}).  As a first observation, proposition \ref{proposition:WeaklyIdempotentCompleteR3} shows that $\EE$ needs to be weakly idempotent complete in order to admit axiom \ref{R3+}.  This indicates that it is not sufficient to just refine the conflation structure on $\EE$ (thus, enlarge the set $\bD$) and that we need to additionally adjust the underlying additive category $\AA$.  A similar remark holds for axiom \ref{R3} (see example \ref{example:ClosingUnderDirectSummandsDoesNotYieldR3}) below.

\begin{definition}\label{definition:R3ClosureUniversalProperty}
	Let $\EE=(\AA,\bD)$ be a deflation-exact category. The \ref{R3}-closure of $\EE$ is a deflation-exact category $\EE_{\textbf{R3}}=(\AA_{\textbf{R3}}, \bD_{\textbf{R3}})$ satisfying axiom \ref{R3} together with an exact functor $\phi\colon \EE\hookrightarrow \EE_{\textbf{R3}}$ satisfying the following 2-universal property: for any strongly deflation-exact category $\FF$, the functor $-\circ \phi\colon \Hom(\EE_{\textbf{R3}}, \FF) \to \Hom(\EE, \FF)$ is an equivalence.
	
	Similarly, one defines the \ref{R3+}-closure.
\end{definition}

The aim of this section is to show the following theorem.

\begin{theorem}
	Let $\EE=(\AA, \bD)$ be a deflation-exact category satisfying axiom \ref{R0*}.
	\begin{enumerate}
		\item The category $\EE$ has an \emph{\ref{R3}-closure} $\EE_{\textbf{R3}}$.
		\item The category $\EE$ has an \emph{\ref{R3+}-closure} $\EE_{\textbf{R3}^+}$.
	\end{enumerate}
	Moreover, the natural maps $\EE\hookrightarrow \EE_{\textbf{R3}}\hookrightarrow \EE_{\textbf{R3}^+}$ lift to triangle equivalences on the bounded derived categories.
	
	Additionally, there exists deflation-exact structure on $\EE$, denoted by $\EE_{\textbf{R3}^-}=(\AA, \bD_{\textbf{R3}^-})$ extending the conflation-structure $\bD$ and satisfying axiom \ref{R3-}. Moreover, $\Db(\EE)\simeq \Db(\EE_{\textbf{R3}^-})$ are triangle equivalent and $\EE_{\textbf{R3}^-}$ is the largest conflation structure on $\EE$ that does not change the bounded derived category up to triangle equivalence.
\end{theorem}

\subsection{The \texorpdfstring{\ref{R3+}}{R3+}-closure}

By proposition \ref{proposition:WeaklyIdempotentCompleteR3}, a deflation-exact category $\EE=(\AA,\bD)$ satisfying axiom \ref{R3+} must be weakly idempotent complete. Given a deflation-exact category $\EE$ satisfying axiom \ref{R0*}, one endows the weak idempotent completion $\widehat{\AA}$ with the following conflation structure $\widehat{\bD}$: a sequence $X \inflation Y \deflation Z$ is a conflation in $\widehat{\AA}$ if and only if it is a direct summand in $\widehat{\AA}$ of a conflation in $\bD$. We simply write $\widehat{\EE}=(\widehat{\AA},\widehat{\bD})$ for the weak idempotent completion of $\EE$ with the above conflation structure. By \cite[appendix~B]{HenrardVanRoosmalen19b}, the following theorem holds:

\begin{theorem}\label{theorem:WeakIdempotentCompletion}
Let $\EE$ be a deflation-exact category satisfying axiom \ref{R0*}. The weak idempotent completion $\widehat{\EE}$ satisfies axiom \ref{R3+} and the exact embedding $\EE \to \widehat{\EE}$ is $2$-universal among exact functors to deflation-exact categories satisfying axiom \ref{R3+}.

Moreover, the natural embedding $\EE \to \widehat{\EE}$ lifts to a triangle equivalence $\Db(\EE) \to \Db(\widehat{\EE})$.
\end{theorem}

It follows that the \ref{R3+}-closure $\EE_{\textbf{R3}^+}$ of a deflation-exact category $\EE$ is simply the weak idempotent completion $\widehat{\EE}$ endowed with the natural conflation structure.

The following corollary states that for a deflation-exact category $\EE$, axiom \ref{R3} implies the dual of theorem \ref{theorem:EquivalentFormulations}.\eqref{item:EquivalentFormulationsD} and axiom \ref{R3+} implies the dual of theorem \ref{theorem:WeaklyIdempotentCompleteEquivalentFormulations}.\eqref{item:WeaklyIdempotentCompleteEquivalentFormulationsF}.

\begin{corollary}\label{corollary:R3GivesWeakL3}
	Let $\EE$ be a deflation-exact category.
	\begin{enumerate}
		\item If $\EE$ satisfies axiom \ref{R3}, then for any morphism $f\colon X\to Y$ having a cokernel $g\colon Y\to Z$ such that there is an inflation $h\colon Y\inflation U$ such that $h\circ f$ is an inflation, $f$ is an inflation.
		\item If $\EE$ satisfies axiom \ref{R3+}, then for any morphism $f\colon X\to Y$ such that there is an inflation $h\colon Y\inflation U$ such that $h\circ f$ is an inflation, $f$ is an inflation.
	\end{enumerate}
\end{corollary}

\begin{proof}
	We show the first statement, the second is similar. One naturally obtains the following commutative diagram:
	\[\xymatrix{
		X\ar[r]^f\ar@{=}[d] & Y\ar[r]^g\ar@{>->}[d]^h & Z\ar@{.>}[d]^k\\
		X\ar@{>->}[r]\ar@{->>}[d] & U\ar@{->>}[r]^m\ar@{->>}[d] & V\ar@{.>}[d]^l\\
		0\ar@{>->}[r] & W\ar@{=}[r] & W
	}\] By theorem \ref{theorem:WeakIdempotentCompletion}, the right column is a conflation in $\widehat{\EE}$ and thus a kernel-cokernel pair in $\EE$. Axiom \ref{R3} implies that the right column is a conflation in $\EE$. Proposition \ref{proposition:MitchellPullback} implies that the upper-right square is a pullback square. By axiom \ref{R2}, the upper row is a conflation.
\end{proof}

The following proposition shows that to any conflation in $\EE_{\textbf{R3}^+}$, one can add split conflations in $\EE$ such that the direct sum is a conflation in $\EE$.

\begin{proposition}\label{proposition:ConflationsOfWeakIdempotentCompletionAreAlmostInOriginalCategory}
	Let $\EE$ be a deflation-exact category satisfying axiom \ref{R0*}. For any conflation $C$ in $\EE_{\textbf{R3}^+}$, there exists a split conflation $C'$ in $\EE$ such that $C\oplus C'$ is a conflation in $\EE$.
\end{proposition}

\begin{proof}
	Let $\overline{X}\inflation \overline{Y}\deflation \overline{Z}$ be a conflation in $\EE_{\textbf{R3}^+}$.  We claim that there exists a conflation $X\inflation Y\deflation Z$ in $\EE$ such that the direct sum $\overline{X}\oplus X\inflation \overline{Y}\oplus Y\deflation \overline{Z}\oplus Z$ is a conflation in $\EE$.
	
Indeed, by definition, there is a conflation $X'\inflation Y'\deflation Z'$ in $\EE_{\textbf{R3}^+}$ such that $\overline{X}\oplus X'\inflation \overline{Y}\oplus Y'\deflation \overline{Z}\oplus Z'$ is a conflation in $\EE$. By remark \ref{remark:ComplementsOfObjectsInWeakIdempotentCompletions}, there are objects $X'_c,Y'_c,Z'_c\in \EE$ such that $X'\oplus X'_c, Y'\oplus Y'_c,Z'\oplus Z'_c\in \EE$. Thus adding the split conflations $X'_c=X'_c\deflation 0, Y'_c=Y'_c\deflation 0, 0\inflation Z'_c=Z'_c$ to the conflation $X'\inflation Y'\deflation Z'$ shows the claim.

By axiom \ref{R2}, taking the pullback of $\overline{Z}\oplus Z\xrightarrow{\begin{psmallmatrix}1&0\\0&0\end{psmallmatrix}} \overline{Z}\oplus Z$ along $\overline{Y}\oplus Y\deflation \overline{Z}\oplus Z$ yields the following commutative diagram in $\EE$:
\[\xymatrix{
	\overline{X}\oplus X\ar@{>->}[r] & \overline{Y}\oplus Y\ar@{->>}[r] & \overline{Z}\oplus Z\\
	\overline{X}\oplus X\ar@{>->}[r]\ar@{=}[u] & \overline{Y}\oplus X\oplus Z\ar@{->>}[r]\ar[u] & \overline{Z}\oplus Z\ar[u]^{\begin{psmallmatrix}1&0\\0&0\end{psmallmatrix}}
}\] The lower row is simply the direct sum (in $\EE_{\textbf{R3}^+}$) of the conflation $\overline{X}\inflation \overline{Y}\deflation \overline{Z}$ and the split conflations $X=X\deflation0$ and $0\inflation Z=Z$. This concludes the proof.
\end{proof}

\subsection{The \texorpdfstring{\ref{R3}}{R3}-closure}

We now turn our attention to the \ref{R3}-closure.  The next example shows that, as is the case for the \ref{R3+}-closure, enlarging the conflation structure alone is, in general, not sufficient.

\begin{example}\label{example:ClosingUnderDirectSummandsDoesNotYieldR3}
Let $Q$ be an equioriented quiver of type $A_3$ and let $\DD = \rep_k Q$ be the category of $k$-linear representations, where $k$ is a field.  The Auslander-Reiten quiver of $\DD$ is given by 
	
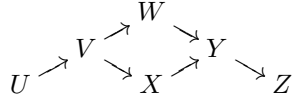
\[\xymatrix@!C=0,5pt@R=1pt{
	&&W\ar[rd]&&\\
	&V\ar[ru]\ar[rd]&&Y\ar[rd]&\\
	U\ar[ru]&&X\ar[ru]&&Z
	}\] Let $\FF$ be the full additive subcategory of $\DD$ of all objects not isomorphic to $V$.  We endow $\FF$ with the structure of a deflation-exact category by taking all conflations $A \inflation B \deflation C$ in $\DD$ where $A$ is not isomorphic to $U.$  It can be explicitly verified that this is indeed a deflation-exact structure (see also remark \ref{remark:ExampleViaQuotient} for a less direct approach). 
	
	The category $\FF$ does not satisfy axiom \ref{R3}. Indeed, there is a conflation $U^{\oplus 2}\inflation W^{\oplus 2}\deflation Y^{\oplus 2}$, and $U\to W\to Y$ is a direct summand.  Note, however, that the pullback of $W \to Y$ along $X \to Y$ does not exist in $\FF$, so that $W \to Y$ cannot be a deflation in any deflation-exact structure on $\FF$.  By theorem \ref{theorem:EquivalentFormulations}, we see that $\FF$ cannot be refined to a deflation-exact structure satisfying axiom \ref{R3}.
\end{example}

Theorem \ref{theorem:EquivalentFormulations}.\eqref{item:EquivalentFormulationsJ} allows a straightforward construction of the \ref{R3}-closure.

\begin{proposition}\label{proposition:ExtensionClosureInWICIsR3Closure}
	Let $\AA_{\textbf{R3}}\subseteq \widehat{\AA}$ be the extension closure of $\AA$ in $\widehat{\EE}=(\widehat{\AA},\widehat{\bD})$. The induced deflation-exact structure $\bD_{\textbf{R3}}$ on $\AA_{\textbf{R3}}$ satisfies axiom \ref{R3} and is the $\ref{R3}$-closure $\EE_{\textbf{R3}}=(\AA_{\textbf{R3}},\bD_{\textbf{R3}})$ of $\EE$.
	
	Moreover, the inclusion $\EE\hookrightarrow \EE_{\textbf{R3}}$ lifts to a triangle equivalence $\Db(\EE)\to \Db(\EE_{\textbf{R3}})$. 
\end{proposition}

\begin{proof}
	By theorem \ref{theorem:EquivalentFormulations}.\eqref{item:EquivalentFormulationsJ}, $\EE_{\textbf{R3}}$ satisfies axiom \ref{R3}. Let $F\colon \EE\to \FF$ be an exact functor to a deflation-exact category $\FF$ satisfying axiom \ref{R3}. The functor $F$ lifts to an exact functor $\widehat{F}\colon \widehat{\EE}\to \widehat{\FF}$. It suffices to show that the restriction of $\widehat{F}$ to $\EE_{\textbf{R3}}$ maps to $\FF\subseteq \widehat{\FF}$. By theorem theorem \ref{theorem:EquivalentFormulations}.\eqref{item:EquivalentFormulationsJ}, $\FF$ lies extension closed in $\widehat{\FF}$, the result follows. 
	
	To see that $\EE\hookrightarrow \EE_{\textbf{R3}}$ lifts to a derived equivalence, it suffices to note that $\EE$ and $\EE_{\textbf{R3}}$ have the same weak idempotent completion. Indeed, the result then follows from theorem \ref{theorem:WeakIdempotentCompletion}.
\end{proof}

It is easy to see that the largest deflation-exact structure $\bDm$ on $\AA$ is given by all semi-stable cokernels (see definition \ref{Definition:SemiStableCokernel}), the next corollary is straightforward.

\begin{corollary}
For a deflation-exact category $\EE = (\AA, \bD)$, we have $\bD_\Rtm = \widehat{\bD} \cap \bDm.$
\end{corollary}

\begin{remark}\label{remark:ExampleViaQuotient}
The deflation-exact category $\FF$ from example \ref{example:ClosingUnderDirectSummandsDoesNotYieldR3} can be obtained as follows. Let $Q$ be the quiver 
	\[\xymatrix@!C=0,5pt@R=1pt{
		1&2\ar[l]_{\alpha}&3\ar[l]_{\beta}&4\ar[l]_{\gamma}
	}\]
	with relation $\alpha\beta=0$.  The category $\rep_k Q$ of finite-dimensional $k$-representations can be visualized via its Auslander-Reiten quiver:
	
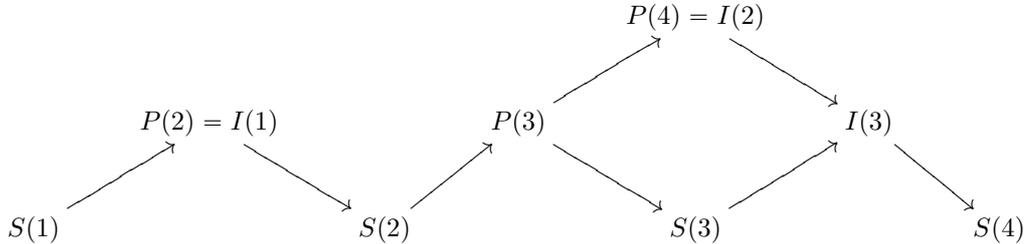
\[\xymatrix{
	 &&&& P(4)=I(2)\ar[rd]&&\\
		& P(2)=I(1)\ar[rd] && P(3)\ar[ru]\ar[rd] && I(3)\ar[rd]&\\
		S(1)\ar[ru] && S(2)\ar[ru] && S(3)\ar[ru] && S(4)\\
	}\]
	Let $\EE$ be the full additive subcategory of $\rep_k Q$ of all objects not isomorphic to $S(1)^{\oplus n}\oplus S(2)$, nor to $S(1)^{\oplus n}\oplus P(3)$ (for any $n \geq 0$).  As $\EE$ is an extension-closed subcategory of $\CC$, it inherits an exact structure. Let $\AA$ be the full additive subcategory of $\EE$ generated by $S(1)$, and let $S_\AA \subset \Mor \EE$ be the set of admissible morphisms with kernel and cokernel in $\AA$.  It is straightforward to check that $\AA\subseteq \EE$ is an admissibly deflation-percolation subcategory in the sense of \cite{HenrardVanRoosmalen19a}, so that the localization $\EE[S_\AA^{-1}]$, equipped with the coarsest conflation structure for which $Q\colon \EE \to \EE[S_\AA^{-1}]$ is exact, is a deflation-exact category.
	
	Note that in $\EE[S_\AA^{-1}]$, we have that $S(1) \cong 0$ and $P(2) \oplus E \cong S(2) \oplus E$, for all nonzero $E \not\in \AA$.  The category $\FF$ is equivalent to $\EE[S_\AA^{-1}]$ and hence a deflation-exact category.
\end{remark}

\subsection{The conflation structure \texorpdfstring{$\EE_{\textbf{R3}^-}$}{R3-}}

Let $\EE$ be a deflation-exact category satisfying axiom \ref{R0*}. In contrast to the \ref{R3+}-closure and the \ref{R3}-closure, for the \ref{R3-}-conflation structure, we do not change the underlying category, but only the conflation structure.  We start with the following definition.

\begin{definition}
	Let $\EE$ be a deflation-exact category satysifying axiom \ref{R0*}. A morphism $f\colon X\to Y$ in $\EE$ is called a \emph{$\mathsf{P}$-deflation} if it satisfies the following two properties:
	\begin{enumerate}[label=\textbf{P\arabic*},start=1]
		\item\label{S1} All pullbacks along $f$ exist.
		\item\label{S2} There exists a map $h$ such that $f\circ h$ is a deflation in $\EE$.
	\end{enumerate}
\end{definition}

\begin{remark}
	\begin{enumerate}
		\item Any deflation is a $\mathsf{P}$-deflation.
		\item If $\EE$ satisfies axiom \ref{R3}, any $\mathsf{P}$-deflation is also a deflation.
	\end{enumerate}
\end{remark}

\begin{lemma}\label{lemma:PullbackOfStrongDeflationIsStrongDeflation}
	$\mathsf{P}$-deflations are stable under pullbacks.
\end{lemma}

\begin{proof}
	This is a straightforward application of the pullback lemma.
\end{proof}

\begin{proposition}\label{proposition:R3MinusClosure}
	Let $\EE$ be a deflation-exact category satisfying axiom \ref{R0*}. The collection of $\mathsf{P}$-deflations in $\EE$ defines a deflation-exact structure on $\EE$ satisfying axiom \ref{R3-}.
\end{proposition}

\begin{proof}
	Let $f\colon X\to Y$ be a $\mathsf{P}$-deflation. By property \ref{S1}, the pullback of $0\to Y$ along $f$ exists and thus $f$ admits a kernel $K$. We first show that $K\to X\stackrel{f}{\rightarrow} Y$ is a kernel-cokernel pair. By property $\ref{S2}$, there is a map $h\colon A\to X$ such that the composition $f\circ h$ is a deflation. Write $K'$ for the kernel of $f\circ h$. Taking the pullback of $f$ along $f\circ h$ it is straightforward to show that we obtain the following commutative diagram:
	\[\xymatrix{
		& K\ar@{=}[r]\ar@{>->}[d] & K\ar[d]\\
		K'\ar@{>->}[r]\ar@{=}[d] & K\oplus A\ar@{->>}[r]\ar@{->>}[d] & X\ar[d]^f\\
		K'\ar@{>->}[r] & A\ar@{->>}[r]_{f\circ h}\ar@{.>}[ru]^h & Y
	}\]Proposition \ref{proposition:MitchellPushout} yields that the lower right square is a pushout as well. It follows that $f$ is the cokernel of $K\to X$ as required.

We now show that the $\mathsf{P}$-deflations define a deflation-exact structure satisfying axiom \ref{R3-}. Clearly $X\to 0$ is a $\mathsf{P}$-deflation as $X\to 0$ is already a deflation by axiom \ref{R0*} in $\EE$. Let $f\colon X\to Y$ and $g\colon Y\to Z$ be two $\mathsf{P}$-deflations. By the pullback lemma, the composition $gf$ satisfies \ref{S1}. By \ref{S2}, there exists a map $h\colon B\to Y$ such that $gh$ is a deflation. By \ref{S1}, the pullback $P$ of $h$ along $f$ exists. Hence we obtain the following commutative diagram:
	\[\xymatrix{
		P\ar[r]^{f'}\ar[d]^{h'} & B\ar@{->>}[r]\ar[d]^h & Z\ar@{=}[d]\\
		X\ar[r]^f & Y\ar[r]^g & Z
	}\]
	By lemma \ref{lemma:PullbackOfStrongDeflationIsStrongDeflation}, the map $f'$ is a $\mathsf{P}$-deflation. Hence there exists a map $h''\colon A\to P$ such that $f'h''$ is a deflation. It follows that $(gf)(h'h'')$ is a deflation by axiom \ref{R1}. Hence $gf$ is a $\mathsf{P}$-deflation. It follows that the $\mathsf{P}$-deflations satisfy axiom \ref{R1}. Lemma \ref{lemma:PullbackOfStrongDeflationIsStrongDeflation} shows that the $\mathsf{P}$-deflations satisfy axiom \ref{R2} as well. By construction, axiom \ref{R3-} is satisfied.
\end{proof}

\begin{corollary}\label{corollary:R3MinusClosure}
Let $\EE =  (\AA, \bD)$ be a deflation-exact category.  There is a smallest deflation-exact structure $\bD_\Rtm (\supseteq \bD)$ on $\AA$ such that $\EE_\Rtm = (\AA, \bD_\Rtm)$ satisfies axiom \ref{R3-}.
\end{corollary}

\begin{proof}
This follows from proposition \ref{proposition:R3MinusClosure} where $\bD_\Rtm$ consists of all conflations whose deflations are $\mathsf{P}$-deflations. 
\end{proof}

The following proposition characterizes the conflation structure $\EE_{\textbf{R3}^-}$ in a universal way.

\begin{proposition}\label{proposition:R3MinusDerivedEquivalence}
	The conflation structure on $\EE_{\textbf{R3}^-}$ is the largest deflation-exact structure on $\EE$ such that the identity $\EE\to\EE_{\textbf{R3}^-}$ is a conflation-exact functor lifting to a triangle equivalence $\Db(\EE)\to \Db(\EE_{\textbf{R3}^-})$.
\end{proposition}

\begin{proof}
	Let $\EE' = (\AA,\bD')$ be a deflation-exact structure on $\AA$ such that $1_\AA\colon \EE\to \EE'$ is an exact functor that lifts to a triangle equivalence $\Db(\EE)\to \Db(\EE')$. It suffices to show that $1_{\EE}\colon \EE'\to \EE_{\textbf{R3}^-}$ is exact. 
	
	Combining theorem \ref{theorem:WeakIdempotentCompletion} with the fact that $1_{\AA}$ lifts to a triangle equivalence $\Db(\EE)\to \Db(\EE')$, we see that the embedding $1_{\AA}\colon\widehat{\EE}\hookrightarrow \widehat{\EE'}$ lifts to a triangle equivalence $\Db(\EE)\stackrel{\sim}{\rightarrow} \Db(\widehat{\EE}')$. Using the universal property of the weak idempotent completion, we obtain the following commutative diagram:
	\[\xymatrix{
		\EE\ar@{^{(}->}[r]\ar@{^{(}->}[d]^{1_{\EE}} & \widehat{\EE}\ar@{^{(}->}[r]\ar@{^{(}->}[d] & \Db(\EE)\ar[d]^{\rotatebox{90}{$\simeq$}}\\
		\EE'\ar@{^{(}->}[r]& \widehat{\EE'}\ar@{^{(}->}[r]& \Db(\widehat{\EE'}) 
	}\] By \cite[proposition~6.2]{HenrardVanRoosmalen19b}, the categories $\widehat{\EE}$ and $\widehat{\EE'}$ have the same conflation structures as conflations correspond to triangles.
	
	Let $X\to Y\to Z$ be a conflation in $\EE'$. As $\widehat{\EE}=\widehat{\EE'}$, the conflation $X\to Y\to Z$ is a retract (in $\widehat{\EE}$) of a conflation $\underline{X}\inflation \underline{Y}\deflation \underline{Z}$ in $\EE$. We obtain the following commutative diagram
	\[\xymatrix{
		X \ar[r]& Y\ar[r] & Z\\
		\underline{X}\ar@{>->}[r]\ar[u] & \underline{Y}\ar@{->>}[r]\ar[u] & \underline{Z}\ar[u]\\
		\underline{X}\ar@{=}[u]\ar@{>->}[r] & P\ar@{->>}[r]\ar[u] & Z\ar[u]\\
		X\ar[r]\ar[u] & Y\ar[u]\ar[r] & Z	\ar@{=}[u]
	}\] where the vertical arrows compose to the identities. Note that we factored the section from $X\to Y\to Z$ to $\underline{X}\inflation \underline{Y}\deflation\underline{Z}$ via the conflation $\underline{X}\inflation P\deflation Y$ obtained by taking the pullback of $Z\to \underline{Z}$ along $\underline{Y}\deflation\underline{Z}$. Note that the composition $P\to \underline{Y}\to Y\to Z$ is a deflation in $\EE$ and thus $Y\to Z$ is a $\mathsf{P}$-deflation. It follows that $1_{\AA}\colon \EE'\to \EE_{\textbf{R3}^-}$ is exact.
\end{proof}

The following example shows that $\EE_{\textbf{R3}^-}$ does not satisfy a 2-universal property as in definition \ref{definition:R3ClosureUniversalProperty}.

\begin{example}
	Let $\FF$ be the deflation-exact category constructed in example \ref{example:ClosingUnderDirectSummandsDoesNotYieldR3}. Let $\EE$ be the deflation-exact subcategory of $\FF$ generated by $U,W$ and $Y$. Note that $W\to Y$ is not a deflation in $\EE$ but it is a deflation in $\EE_{\textbf{R3}^-}$. It follows that the exact embedding $\EE\hookrightarrow \FF$ does not lift to an exact map $\EE_{\textbf{R3}^-}\hookrightarrow \FF_{\textbf{R3}^-}$.
\end{example}

\begin{remark}
	\begin{enumerate}
		\item The \ref{R3-}-closure of a (possibly non-additive) deflation-exact category has been introduced in \cite[\S1.4]{Rosenberg11} as the \emph{closure} of the category $\EE$.
		\item	Proposition \ref{proposition:R3MinusDerivedEquivalence} extends the fact that the \ref{R3-}-closure of a deflation-exact category preserves projectives (see \cite[proposition~1.4.4]{Rosenberg11}).
	\end{enumerate}
\end{remark}

\subsection{About the conflations in the obscure closures}

We now have a closer look at which sequences $X \to Y \to Z$ in a deflation-exact category $\EE$ become conflations in the closures considered in this section.  The following proposition is the main result.

\begin{proposition}\label{proposition:WhenConflationInCompletion}
Let $X\xrightarrow{f} Y \xrightarrow{g} Z$ be a sequence in a deflation-exact category $\EE$.  The following are equivalent:
\begin{enumerate}
	\item\label{enumerate:ConflationInClosure1} the sequence is a conflation in $\EE_\Rt$,
	\item\label{enumerate:ConflationInClosure2} the sequence is a conflation in $\EE_\Rtp$,
	\item\label{enumerate:ConflationInClosure3} the sequence is a conflation in $\ex{\EE}$,
	\item\label{enumerate:ConflationInClosure4} there is an $A \in \EE$ such that $\xymatrix@1{X\oplus A\ar@{>->}[r]^-{\begin{psmallmatrix} f & 0 \\ 0 & 1 \end{psmallmatrix}} & Y \oplus A \ar@{->>}[r]^-{\begin{psmallmatrix} g & 0 \\ \end{psmallmatrix}} & Z}$ is a conflation in $\EE.$
\end{enumerate}
If $g\colon Y \to Z$ admits all pullbacks, then the previous are furthermore equivalent to
\begin{enumerate}[resume]
  \item\label{enumerate:ConflationInClosure5} the sequence is a conflation in $\EE_\Rtm$.
\end{enumerate}
\end{proposition}

\begin{proof}
If \eqref{enumerate:ConflationInClosure4} holds, then \eqref{enumerate:ConflationInClosure2} holds, as $X \to Y\to Z$ is a direct  summand of a  conflation in $\EE$. If \eqref{enumerate:ConflationInClosure2} holds, then there is a sequence $A \to B \to C$ such that $X \oplus A \to Y \oplus B \to Z \oplus C$ is a conflation in $\EE$.  The pullback of  this conflation along the embedding $Z \to Z \oplus C$ yields the required conflation $X \oplus A \inflation Y \oplus A \deflation Z,$ and hence \eqref{enumerate:ConflationInClosure4} holds.

The equivalence $\eqref{enumerate:ConflationInClosure1} \Leftrightarrow \eqref{enumerate:ConflationInClosure2}$ holds as $\EE_\Rt$ lies extension-closed in $\EE_\Rtp.$  The equivalence $\eqref{enumerate:ConflationInClosure1} \Leftrightarrow \eqref{enumerate:ConflationInClosure3}$ is established in theorem \ref{theorem:EquivalentFormulations}.

Finally, assume that $g\colon Y \to Z$ is a semi-stable cokernel.  If \eqref{enumerate:ConflationInClosure4} holds, then it follows from proposition \ref{Proposition:ObscureAxiomsCharacterizationsWithSummands} that \eqref{enumerate:ConflationInClosure5} holds.  The other direction follows since conflations are closed under direct sums.
\end{proof}

As a corollary, we obtain the following description of the closure $\EE_\Rtm.$

\begin{corollary}\label{corollary:R3mClosureViaR3p}
For a deflation-exact category $\EE = (\AA, \bD)$, we have $\bD_{\Rtm} = \bDm \cap \bD_\Rtp.$
\end{corollary}

\begin{proposition}
Let $\AA$ be an additive category with two deflation-exact structures $\bD$ and $\bD'$.  The following are equivalent:
\begin{enumerate}
	\item\label{enumerate:WhenEqualClosures1} $(\AA, \bD)_{\Rtm} = (\AA, \bD')_{\Rtm}$,
	\item\label{enumerate:WhenEqualClosures2} $(\AA, \bD)_{\Rt} = (\AA, \bD')_{\Rt}$,
	\item\label{enumerate:WhenEqualClosures3} $(\AA, \bD)_{\Rtp} = (\AA, \bD')_{\Rtp}$.
\end{enumerate}
\end{proposition}

\begin{proof}
The implication $\eqref{enumerate:WhenEqualClosures3} \Rightarrow \eqref{enumerate:WhenEqualClosures2}$ holds as $\EE_\Rt$ is the extension-closure of $\EE$ in $\EE_\Rtp.$  For the implication $\eqref{enumerate:WhenEqualClosures2} \Rightarrow \eqref{enumerate:WhenEqualClosures1}$, it follows from proposition \ref{proposition:WhenConflationInCompletion} that a sequence $X \to Y \xrightarrow{g} Z$ in $\AA$ lies in $(\AA, \bD)_{\Rtm}$ if and only if it lies in $(\AA, \bD)_{\Rt}$ and $g$ is a semi-stable cokernel.  Finally, for the implication $\eqref{enumerate:WhenEqualClosures1} \Rightarrow \eqref{enumerate:WhenEqualClosures3}$, it suffices to notice that $(\EE_{\Rtm})_\Rtp = \EE_\Rtp.$
\end{proof}

We end this section with a variant of proposition \ref{proposition:WhenConflationInCompletion}.

\begin{proposition}\label{proposition:WhenConflationInWic}
Let $\EE$ be a deflation-exact category.  A sequence $X \to Y \to Z$ in $\EE_{\Rtp}$ is a conflation if and only if there is a split conflation $A \to B \to C$ in $\EE$ such that $A \oplus X \inflation B \oplus Y \deflation C \oplus Z$ is a conflation in $\EE.$
\end{proposition}

\begin{proof}
Let $X \inflation Y \deflation Z$ be a conflation in $\EE_{\Rtp}$.  By the definition of the weak idempotent completion, there is an object $X' \in \EE$ such that $X \oplus X' \in \EE$.  Similar complements $Y', Z' \in \EE$ exist for $Y,Z \in \EE.$  We find that the conflation
\[(X'\oplus Y') \oplus X \inflation (X'\oplus Y' \oplus Z') \oplus Y \deflation Z' \oplus Z\]
in $\EE_{\Rtp}$  of  which  all terms lie in $\EE.$  The existence of the split sequence $A \to B \to C$ as in the statement of the proposition now follows from proposition \ref{proposition:WhenConflationInCompletion}.

The converse follows from the definition of the the conflation structure of $\EE_{\Rtp}$.
\end{proof}
\section{The lattice of strongly one-sided exact structures}\label{Section:Lattices}

Let $\AA$ be an additive category.  It is known that the set of exact structures on an additive category $\AA$ forms a complete lattice (see for example \cite{BrustleHassounLangfordRoy20, FangGorsky20}), and it is easy to see that the same observation holds for one-sided exact structures (see \cite{Rosenberg11}).  These observations have been generalized in \cite{BaillargeonBrustleGorskyHassoun20} to the setting of weakly (inflation- or deflation-)exact structures.

In this section, we have a closer look at the posets of deflation-exact structures satisfying the various versions of the obscure axiom from definition \ref{definition:VariantsOfObscureAxioms}.  We show that the lattices of deflation-exact structures satisfying axiom \ref{R3} or \ref{R3-} are ideals in the strong deflation-exact structures on the weak idempotent completion $\widehat{\AA}.$ 

We start by recalling the notion of a lower and an upper adjoint of a morphism between posets.

\begin{definition}
	Let $(A,\leq)$ and $(B,\leq)$ be two partially ordered sets. A \emph{(monotone) Galois connection} between $(A,\leq)$ and $(B,\leq)$ consists of two monotone functions $F\colon A\to B$ and $G\colon B\to A$ such that for all $a\in A$ and $b\in B$, we have 
	\[F(a)\leq b \Leftrightarrow a\leq G(b).\]
	The map $F$ is called the \emph{lower adjoint} and the map $G$ is called the \emph{upper adjoint}.  We write $F = G^\flat$ and $G = F^\sharp.$
\end{definition}

The following lattices will be considered in this section.

\begin{notation}
	Let $\AA$ be a (small) additive category.  We write $\fD(\AA)$ for set of deflation-exact structures on $\AA$.  The set of strongly deflation-exact structures on $\AA$ is denoted by $\sfD(A)$.  We write $\fDm(\AA)$ and $\fDp(\AA)$ for the sets of deflation-exact structures satisfying axiom \ref{R3-} or \ref{R3+}, respectively.
	
	Likewise, the set of inflation-exact structures on $\AA$ is denoted by $\fI(\AA)$ and the set of strongly inflation-exact structures on $\AA$ is denoted by $\mathfrak{sI}_{\AA}$.
	
	The partially ordered set of (two-sided) exact structures on $\AA$ is denoted by $\fE(\AA)$.
\end{notation}

\begin{proposition}
Let $\AA$ be an additive category.  The posets $\fD(\AA), \sfD(\AA),$ and $\fDm(\AA)$ (ordered by inclusion) are complete lattices.  If $\AA$ is weakly idempotent complete (so that $\fDp(\AA)$ is nonempty), then $\fDp\AA)$ is a complete lattice.
\end{proposition}

\begin{proof}
The meet operation is given by the intersection.  The least element is given by the split conflation structure (this satisfies axiom \ref{R3+} if and only if $\AA$ is weakly idempotent complete). The greatest element of $\fD$ is given by \cite[proposition 1]{Rump11} (this is also the greatest element of $\fDm$ by corollary \ref{corollary:R3MinusClosure}), and the greatest element of $\sfD$ is given by \cite[corollary 1]{Rump11}.  When $\AA$ is weakly idempotent complete, a greatest element of $\fDp(\AA)$ is given by \cite{Crivei12}.
\end{proof}

\begin{proposition}\label{proposition:ObscureAxiomsInherited}
Let $\bD \subseteq \bD'$ be two deflation-exact structures on an additive category $\AA$.  Assume that $\bD$ satisfies axiom \ref{R3-}.
\begin{enumerate}
	\item If $\bD'$ satisfies axiom \ref{R3} or \ref{R3+}, then so does $\bD.$
	\item If $\bD'$ is exact, then so is $\bD.$
\end{enumerate}
\end{proposition}

\begin{proof}
The first statement is straightforward to verify.  For the second statement, note that $\bD$ satisfies \ref{R3} as $\bD'$ does.  It then follows from \cite[theorem~1]{Rump11} that $\bD = \bD \cap \bD'$ is an exact structure on $\AA.$
\end{proof}

\begin{remark}
The previous proposition can be reformulated as follows: the posets $\sfD(\AA)$ and $\fE(\AA)$ are ideals in $\fDm(\AA).$
\end{remark}

\begin{notation}
Let $\bC$ be any conflation structure on an additive category $\AA$.  We write $\widehat{\bC}$ for the conflation structure on the weak idempotent completion $\widehat{\AA}$ given by: a sequence $X \to Y \to Z$ in $\widehat{\AA}$ lies in $\widehat{\bC}$ if and only if it is a direct summand of a conflation in $\bC$.

We write $\widehat{\fDm(\AA)}$ for the image of $\fDm(\AA)$ under $\bC \mapsto \widehat{\bC}$.  The sets $\widehat{\sfD(\AA)}$ and $\widehat{\fE(\AA)}$ are defined in a similar fashion.
\end{notation}

\begin{remark}
For a deflation-exact category $\EE = (\AA, \bC)$  satisfying axiom \ref{R0*}, we have $\EE_\Rtp = (\widehat{\AA}, \widehat{\bC}).$
\end{remark}

\begin{corollary}
	Let $\AA$ be an additive category. The map $\fDm(\AA) \to \fDp(\widehat{\AA})\colon \bD \mapsto \widehat{\bD}$ is an injection. In particular,	the sets $\widehat{\fDm(\AA)}$, $\widehat{\sfD(\AA)}$, and $\widehat{\fE(\AA)}$ are ideals in $\fDp(\widehat{\AA})$.  Each of the embeddings of the ideals into $\fDp(\widehat{\AA})$ has an upper adjoint.
\end{corollary}

\begin{proof}
	The map $\fDm(\AA) \to \fDp(\widehat{\AA})\colon \bD \mapsto \widehat{\bD}$ is an injection by corollary \ref{corollary:R3mClosureViaR3p}.  This implies that both $\sfD(\AA)$ and $\fE(\AA)$ embed in $\fDp(\widehat{\AA})$ as well.  Each of the embeddings into $\fDp(\widehat{\AA})$ has an upper adjoint, given by intersecting with the appropriate maximal structure on $\AA$ (this follows from proposition \ref{proposition:ObscureAxiomsInherited}).
\end{proof}

\begin{remark}
It is shown in \cite{Rump11} that the embedding $\alpha\colon \sfD(\AA) \to \fD(\AA)$ has an upper adjoint, given by the operator $PQ$ (explicitly, every deflation-exact structure $\bD$ on $\AA$ admits a largest deflation-exact substructure satisfying axiom \ref{R3}).  The embedding does not need to have a lower adjoint, for example, the maximal exact structure from example \ref{example:AboutLatticeOfObscureStructures} need not be embedded into a strong deflation-exact structure.

It follows from corollary \ref{corollary:R3MinusClosure} that the embedding $\beta\colon \fDm(\AA) \to \fD(\AA)$ admits a lower adjoint, mapping $(\AA, \bD)$ to $(\AA, \bD)_{\Rtm}$.  Moreover, $\alpha$ admits a lower adjoint (so every deflation-exact structure $\bD$ on $\AA$ has an obscure closure in $\AA$) if and only if the maximal deflation-exact structure $\bDm$ satisfies axiom \ref{R3}, thus when $\fDm(\AA) = \sfD(\AA).$
\end{remark}

The following example shows that, for an additive category $\AA$ with the maximal deflation-exact structure $\bDm$, the deflation-exact structure $\widehat{\bDm}$ need not be the maximal exact structure on the weak idempotent completion $\widehat{\AA}.$ 

\begin{example}\label{example:AboutLatticeOfObscureStructures}
Let $R = \bC[[ t ]]$ be the formal power series ring in one variable and let $\AA$ be the additive subcategory of $\modd R$ generated by the indecomposables $R$ and $R/(t^n)^{\oplus m}$, for each $m \geq n \geq 1$.  Here, $\widehat{\AA} \simeq \modd R$, so the maximal deflation-exact structure on $\widehat{\AA}$ is given by all kernel-cokernel pairs.

We claim that $R \xrightarrow{t} R \to R/(t)$ is not a conflation in $\widehat{\bDm}$ so that $\widehat{\bDm}$ is  not the maximal deflation-exact structure on $\widehat{\AA}.$  Indeed, for $R \xrightarrow{t} R \to R/(t)$ to be a conflation in $\widehat{\bDm},$ there needs to be a (split) sequence $A \to B \to C$ such that $R \oplus A \inflation R \oplus B \deflation R/(t) \oplus C$ is a conflation in $\bDm$, see proposition \ref{proposition:WhenConflationInWic}.  Let $p\colon R/(t^n) \to R/(t)$ be a nonzero morphism and consider the following commutative diagram
\[\xymatrix{
R \oplus A \ar@{>->}[r] \ar@{=}[d] & R \oplus A \oplus R/(t^{n-1}) \oplus R/(t^n)^{n-1} \ar@{->>}[r] \ar[d] & R/(t^n) \oplus R/(t^n)^{n-1} \ar[d]^{\begin{psmallmatrix} p & 0\\ 0 & 0 \end{psmallmatrix}} \\
R \oplus A \ar@{>->}[r] & R \oplus B \ar@{->>}[r] & R/(t) \oplus C
}\]
where the rows are conflations and the right square is a pullback.  Note that $R \oplus A \oplus R/(t^{n-1}) \oplus R/(t^n)^{n-1}$ can only lie in $\AA$ if $R/(t^n)$ is a direct summand of $A.$  As this can only be the case for finitely many choices of $n$, this shows that not all pullbacks lie in $\AA.$  Consequently, $R \xrightarrow{t} R \to R/(t)$ is not a conflation in $\widehat{\bDm}$ and hence $\widehat{\bDm}$ is not the maximal exact structure on $\widehat{\AA} \simeq \modd R.$
\end{example}

\section{The exact structure of stable kernel-cokernel pairs}\label{Section:StableConflations}

Let $\AA$ be an additive category.  The set $\bC$ of stable kernel-cokernel pairs in $\AA$ is given by the intersection of the maximal inflation-exact and deflation-exact structures on $\AA.$  It is shown in \cite{Rump15} that $\bC$ need not be an exact structure on $\AA$.  Indeed, the maximal exact structure on $\AA$ can be considerably smaller (see \cite{Rump15}). In this subsection, rather than restricting $\bC$ to an exact structure on $\AA$, as is done in \cite{Rump11}, we enlarge the category $\AA$ so that $\bC$ generates an exact structure, that is, we construct an exact category $(\AA', \bE)$ which is, in some way, the smallest exact category containing the conflation category $(\AA, \bC).$  We refer to proposition \ref{proposition:UniversalPropertyIntersection} below for a precise statement.

We proceed in slightly more generality.  Instead of the starting from the intersection of the maximal inflation-exact and deflation-exact structures, we allow to begin with the intersection of any inflation-exact structure $\bI$ and any deflation-exact structure $\bD$.  The category $\AA'$ we construct is a subcategory of the weak idempotent completion $\widehat{\AA}$ and is dependent on the choices of $\bI$ and $\bD.$

\begin{definition}\label{Definition:SemiStableCokernel}
Let $\AA$ be an additive category.  We say that a cokernel $g\colon Y \to Z$ is semi-stable if for every morphism $s\colon Y' \to Z$, the pullback
\[\xymatrix{P \ar@{..>}[r]^p \ar@{..>}[d] & Y' \ar[d]^s \\ X \ar[r]^g & Z}\]
exists and $p\colon P \to Y'$ is itself a cokernel.  A semi-stable kernel is defined dually.

A kernel-cokernel pair $X \xrightarrow{f} Y \xrightarrow{g} Z$ is called stable if $f$ is a semi-stable kernel and $g$ is a semi-stable cokernel.
\end{definition}

The following is shown in \cite{Rosenberg11} and \cite{Rump11}.

\begin{proposition}
Let $\AA$ be an additive category.  The conflation structure $\bC$ consisting of all kernel-cokernel pairs $X \xrightarrow{f} Y \xrightarrow{g} Z$ where $g$ is a semi-stable cokernel is a deflation-exact structure on $\AA$ and is maximal with that property.
\end{proposition}

We continue by having a closer look at $\widehat{\bD} \cap \widehat{\bI}.$

\begin{lemma}\label{Lemma:IntersectionLeftRight}
Let $\bD$ and $\bI$ be a deflation- and an inflation-exact structure on an additive category $\AA$, respectively.  We have
\[\widehat{\bD} \cap \widehat{\bI} = \widehat{\bD \cap \bI}.\]
\end{lemma}

\begin{proof}
The inclusion $\widehat{\bD \cap \bI} \subseteq \widehat{\bD} \cap \widehat{\bI}$ is clear.  For the other inclusion, let $X \inflation Y \deflation Z$ be a conflation in $\widehat{\bD} \cap \widehat{\bI}$.  It follows from proposition \ref{proposition:WhenConflationInWic} that there is a split kernel-cokernel pair $A \inflation B \deflation C$ in $\AA$ such that $X \oplus A \inflation Y \oplus B \deflation Z  \oplus C$ is a conflation in $\bD$, as well as in $\widehat{\bI}$.  The dual of proposition \ref{proposition:WhenConflationInWic} yields the existence of a split kernel-cokernel pair $A' \inflation B' \deflation C'$ such that the conflation $X \oplus (A\oplus A') \inflation Y \oplus (B \oplus B')\deflation Z \oplus (C\oplus C')$ is a conflation in $\bD \cap \bI.$
\end{proof}

\begin{proposition}\label{proposition:IntersectionIsOnWicIsExact}
Let $\bD$ and $\bI$ be a deflation- and an inflation-exact structure on an additive category $\AA.$  The conflation structure $\widehat{\bD} \cap \widehat{\bI} =\widehat{\bD \cap \bI}$ on the weak idempotent completion $\widehat{\AA}$ is an exact structure.
\end{proposition}

\begin{proof}
By theorem \ref{theorem:WeaklyIdempotentCompleteEquivalentFormulations}, the deflation-structure $\widehat{\bD}$ satisfies axiom \ref{R3+}, and dually $\widehat{\bI}$ satisfies axiom \ref{L3+}.  As $\widehat{\bD} \cap \widehat{\bI} = \widehat{\bD \cap \bI}$, it follows from \cite[theorem 1]{Rump11} that $\widehat{\bD} \cap \widehat{\bI}$ is an exact structure on $\widehat{\AA}.$
\end{proof}

\begin{proposition}\label{proposition:UniversalPropertyIntersection}
Let $\bD$ and $\bI$ be a deflation- and an inflation-exact structure on an additive category $\AA.$  There is an exact category $(\AA', \bE)$ and a conflation-exact functor $\varphi\colon (\AA, \bD \cap \bI) \to (\AA', \bE)$ satisfying the following 2-universal property: for each exact category $(\BB, \EE)$, the functor
\[-\circ \varphi\colon \Homex((\AA,\bD \cap \bI),(\BB,\bF)) \to \Homex((\AA',\bE), (\BB,\bF))\]
is an equivalence. 
\end{proposition}

\begin{proof}
The 2-universal property of the weak idempotent completion gives an equivalence $\Hom(\AA,\widehat{\BB}) \to \Hom(\widehat{\AA},\widehat{\BB}).$  Let $\AA'$ be the extension-closure of $\AA$ in $\widehat{\AA}$ with respect to the conflation structure $\widehat{\bD \cap \bI}.$  It follows from proposition \ref{proposition:IntersectionIsOnWicIsExact} that $\widehat{\bD} \cap \widehat{\bI}$ is an exact structure on $\widehat{\AA}$. As $\AA'\subseteq \widehat{\AA}$ is extension-closed, $\AA'$ inherits an exact structure $\bE$ from $(\widehat{\AA},\widehat{\bD} \cap \widehat{\bI})$. To verify that the universal property is satisfied, we consider the following (essentially) commutative diagram
\[\xymatrix{
{(\AA, \bD \cap \bI)} \ar[r] \ar[d] & {(\BB, \bF)} \ar[d] \\
{(\widehat{\AA}, \widehat{\bD} \cap \widehat{\bI})} \ar[r] & {(\widehat{\BB}, \widehat{\bF})} }\]
We only need to verify that the functor $\AA' \to \widehat{\AA} \to \widehat{\BB}$ factors through $\BB$.  However, the essential image of this functor is $\BB \subseteq \widehat{\BB}$ (this uses that $(\BB, \bF)$ lies extension-closed in $(\widehat{\BB}, \widehat{\bF})$, see theorem \ref{theorem:EquivalentFormulations}).  The required property follows from this.
\end{proof}

\begin{remark}
In the proof of proposition \ref{proposition:UniversalPropertyIntersection}, we construct the exact category $(\AA', \bE)$ as the extension-closure of $\AA$ in the intersection of the categories $(\AA, \bI)_{\Ltp}$ and $(\AA, \bD)_{\Rtp}.$  Similarly, one can obtain the exact category $(\AA', \bE)$ as the intersection of $(\AA, \bI)_{\Lt}$ and $(\AA, \bD)_{\Rt}$, seen as subcategories of $\widehat{\AA}$.
\end{remark}

\begin{remark}\label{Remark:InteresectionIsDerivedEquivalent}
As $(\AA, {\bD} \cap {\bI})$ is just a conflation category, its derived category is not defined. Therefore, we cannot express that the closure in proposition \ref{proposition:UniversalPropertyIntersection} gives a derived equivalence.  However, it is easy to see that $\Db(\AA', \bE) = \Kb(\AA) / [\Ac(\AA, \bD) \cap \Ac(\AA, \bI)]$.
\end{remark}

\begin{example}
We return to the setting of example \ref{example:AboutLatticeOfObscureStructures}.  It follows from proposition \ref{proposition:WhenConflationInWic} that, for every $n \geq 2$, the sequence $R/(t)\to R/(t^n) \to R/(t^{n-1})$ is a conflation in $\EE_\Rtp$ (with $A = R/(t^n)^{\oplus n}$ and $C = R/(t^{n-1})^{\oplus (n-1)}$).  Likewise, considering the maximal inflation-exact structure on $\AA$ we find that the sequence $R/(t)\to R/(t^n) \to R/(t^{n-1})$ is a conflation in $\widehat{\bIm}$.   This implies that $\AA' = \widehat{\AA}.$  As the sequence $R \xrightarrow{t} R \to R/(t)$ is not a conflation in $\widehat{\bDm}$, we see that $(\AA', \bE')$ is not the maximal exact structure on $\widehat{\AA} = \modd R.$
\end{example}


\providecommand{\bysame}{\leavevmode\hbox to3em{\hrulefill}\thinspace}
\providecommand{\MR}{\relax\ifhmode\unskip\space\fi MR }
\providecommand{\MRhref}[2]{%
  \href{http://www.ams.org/mathscinet-getitem?mr=#1}{#2}
}
\providecommand{\href}[2]{#2}

\end{document}